\documentclass{amsart}
\usepackage{macros}

\begin{document}
\title[Invertibility of sheaf-theoretic line operators]{Towards invertibility of derived categories of sheaf-theoretic line operators}

\author{Claudia Scheimbauer}
\author{Jackson Van Dyke}
\address{TU Munich, Boltzmannstr. 3, 85746 Garching}
\email{scheimbauer@tum.de}
\email{jackson.van-dyke@tum.de}

\date{\today}

\begin{abstract}
We investigate when the braided monoidal $\infty$-categories of sheaf-theoretic line operators of $3$-dimensional topological quantum field theories are invertible
in the higher Morita category, so that they determine invertible $4$-dimensional anomaly theories.
We consider a zoo of different examples which are inherently derived.
First, enhancing Ben-Zvi--Francis--Nadler, we show that a symmetric monoidal category of quasi-coherent sheaves $\QC\left(X\right)$ is invertible only if $X$ is a point, and convolution categories on loop stacks are Morita trivial whenever the stack satisfies descent along its diagonal, e.g.\ for classifying stacks.
Secondly, we provide an example showing that, unlike for finite braided tensor categories,
non-degeneracy does not imply invertibility in the derived setting: For an odd vector space
$V\left[1\right]$, the $2$-shifted $\cD$-modules $\DMod^2\left(V\left[1\right]\right)$ are $3$-dualizable and non-degenerate
but neither factorizable nor cofactorizable.
This category is closely related to the line operators of linear Rozansky--Witten theory.
Finally, we define an anomaly of a group action on
an invertible braided monoidal $\infty$-category as a Witt class over $\Rep\left(G\right)$, and show that the category of representations of the quantum group of
a torus at a root of unity is invertible relative to its quantum Frobenius.
\end{abstract}

\maketitle

\tableofcontents

\section{Introduction}
\label{sec:intro}

Derived categories of line operators of $3$-dimensional topological quantum field theories (TQFTs) form braided monoidal, i.e.~$\mathbb{E}_2$-monoidal $\infty$-categories, and often arise via
deformation quantization of $2$-shifted symplectic stacks.
As such, they are candidates for determining $4$-dimensional theories.
More explicitly, if an extended $3$-dimensional theory assigns some object $B$ to the point,
then its line operators are its value on the circle, which is the $\EE_2$-algebra
$\mathcal{Z}_1\left(B\right)$.
The object $B$ itself can be viewed relative to $\mathcal{Z}_1\left(B\right)$, hence we can ask when it determines a boundary theory for the
associated $4$-dimensional theory.

In this article we consider a zoo of different examples of $B$ which are suitable sheaf-theoretic categories and as such inherently derived.
We investigate when $\mathcal{Z}_1\left(B\right)$ is {\em invertible}, so that
the $4$-dimensional theory is invertible and is an anomaly of the $3$-dimensional theory. Concretely, $\mathcal{Z}_1\left(B\right)$ is an object in the higher Morita category $\umor_2\left(\Pre\right)$, and we discuss invertibility as such. By \cite{BJSS, SSS} (Theorem \ref{thm:inv_E2}), this holds if and only if it is $3$-dualizable,
non-degenerate, factorizable and cofactorizable; these conditions are phrased using factorization homology.

The examples are motivated by Rozansky--Witten theory, which takes as input a holomorphic
symplectic space, or more generally a $2$-shifted symplectic stack.
Such a space is locally modeled on a symplectic vector space, so the global theory should
be determined by the linear theory and the action of its symmetries (providing
the transition data).
We expect this action to be anomalous, with the anomaly obstructing the gluing.
This is why we focus on the seemingly simple linear case, and why invertibility is the
question of interest: As explained above, the line operators define an anomaly precisely
when they are invertible.
Quasi-coherent sheaves serve as a warm-up for shifted $\cD$-modules, which provide an avatar
for the line operators of linear Rozansky--Witten theory.
The same framework includes Chern--Simons theory, whose input is the $2$-shifted symplectic
stack $BG$.
When $G$ is a torus, this is a global example which is understood, also in the derived
setting, through its quantum group.

\subsection*{Quasi-coherent sheaves}
As a warm-up, we revisit quasi-coherent sheaves and the formalism of Ben-Zvi, Francis and Nadler in \cite{BZFN,BZFN:app}.
Using their formulas for factorization homology we rephrase the invertibility conditions geometrically.
It follows immediately that a commutative algebra is invertible only if its unit is an
equivalence, so that $\left(\QC\left(X\right),\otimes\right)$ is invertible only if $X$ is a point.
Instead, a basic source of invertible
examples is convolution categories on loop stacks $\left(\QC\left(\cL X\right),*\right)$, which are the Drinfeld centers of $\left(\QC\left(X\right),\otimes\right)$.
In order to calculate their factorization homologies, we pass through 
a Morita equivalence between 
sheaves on the loop stack $\cL X= X\times_{X\times X} X$ with convolution and sheaves on the base $X\times X$ with the symmetric monoidal tensor product. 
For loop stacks, under the hypothesis of diagonal descent, we upgrade the Morita equivalence of \cite{BZFN:app} to $\umor_2(\Pre)$, using the criterion of \cite[Theorem 2.26]{BJSS}.

\begin{introthm}[Proposition \ref{prop:inv_symm}, Theorems \ref{thm:inv_rel} and \ref{thm:trivial_QC}]
\label{thm:intro_A}
\leavevmode
\begin{enumerate}
\item For $d\geq1$ and $\coef$ as in Section \ref{sec:mor}, a commutative algebra in $\coef$ is
invertible in $\umor_d\left(\coef\right)$ if and only if its unit is an equivalence. In particular, if $X$ is a quasi-compact stack with affine
diagonal and a $k$-point, then $\left(\QC\left(X\right),\otimes\right)$ is invertible in
$\umor_d\left(\Pre\right)$ only if $X\simeq\pt$.
\item If $X$ satisfies diagonal descent, e.g.\ $X=BG$ or the formal completion $\widehat V_0$ of a
vector space at the origin, then $\left(\QC\left(\cL X\right),*\right)$ is Morita trivial,
trivialized by $\QC\left(X\right)$.
\end{enumerate}
\end{introthm}

Part (2) is the derived-geometric analogue of the Witt triviality of Drinfeld centers of fusion
categories \cite{DMNO}. Diagonal descent is restrictive, however: It fails for any
positive-dimensional scheme (Remark \ref{rmk:diag_descent_fails}), and in general we do not
expect $\QC\left(\cL X\right)$ to be invertible.

\subsection*{Shifted D-modules}
Next, we turn to our original motivation for this article, namely studying the line operators of Rozansky--Witten theory.
Choosing the target $T^*\left[2\right]X$, these are closely related to the
deformation quantization $\DMod^2\left(X\right)$ of the $2$-shifted symplectic stack
$T^*\left[2\right]X$ \cite{Ber:Sph, HoLi} (Section \ref{sec:dmod}).
We show that it is typically not invertible.
To test the invertibility conditions, we compute the higher centers of $\DMod^2\left(X\right)$ using the
excision theorem of \cite{HoLi}, together with rigidity of $\DMod^n\left(X\right)$ for $X$ bounded (Proposition \ref{prop:rigid}), which is observed for $n=1$ in
\cite{Ber:center} and used in \cite[Section 5.2]{Ber:Sph}. Rigidity also gives $3$-dualizability of $\DMod^2\left(X\right)$ for $X$
smooth (Corollary \ref{cor:3dual}).
The symmetric monoidal category $\DMod\left(X\right)$ of ordinary $\cD$-modules is central in
$\DMod^2\left(X\right)$, so $\DMod^2\left(X\right)$ sheafifies over the de Rham stack of $X$
(Remark \ref{rmk:DMod_center}).
More precisely, for $X=BG$, $\DMod^2\left(X\right)$ is non-degenerate exactly when $G$ is
unipotent (Theorem \ref{thm:centers} and Proposition \ref{prop:BG}). Even then, as for an odd vector
space $X=V\left[1\right]$, invertibility can still fail through (co)factorizability.
We therefore expect $\DMod^2\left(X\right)$ to be invertible only when $X$ is a point.

\begin{introthm}[Theorems \ref{thm:centers} and \ref{thm:BW}, Proposition \ref{prop:BG} and Corollary \ref{cor:BG_not_inv}]
\label{thm:intro_B}
\leavevmode
\begin{enumerate}
\item For a smooth stack $X$ as in Section \ref{sec:dmod} and $j=0,1,2$ (with $X\to\cL X$ bounded if $j=2$, e.g.\ $X=BG$), the $\EE_j$-center of $\DMod^2\left(X\right)$ is its
factorization homology over the $j$-sphere:
\begin{equation*}
\cZ_j\left(\DMod^2\left(X\right)\right)\simeq\int_{S^j\times\RR^{2-j}}\DMod^2\left(X\right)
\simeq\DMod^{2-j}\left(X^{S^j}\right) \ .
\end{equation*}
\item\label{thm:intro_B-example} For a non-zero vector space $V$, $\DMod^2\left(V\left[1\right]\right)$ is $3$-dualizable and
non-degenerate, but neither factorizable nor cofactorizable. In particular, it is not invertible in
$\umor_2\left(\Pre\right)$.
\item If $G$ is a non-trivial affine algebraic group which is either not unipotent or abelian, then
$\DMod^2\left(BG\right)$ is not invertible.
\end{enumerate}
\end{introthm}

For finite braided tensor categories over an algebraically closed field of characteristic zero,
non-degeneracy, factorizability and cofactorizability are equivalent to each other and to invertibility \cite{BJSS}. By Theorem \ref{thm:intro_B}(\ref{thm:intro_B-example}), in the
derived setting non-degeneracy implies neither of the other two, even for $3$-dualizable
$\EE_2$-algebras. The reason is visible already in Theorem \ref{thm:intro_B}(1): For $X=V\left[1\right]$,
the mapping stacks $X^{S^0}$ and $X^{S^2}$ are purely odd, whereas
$X^{S^1}\simeq V\times V\left[1\right]$ has an even direction.

\subsection*{Equivariantization}
Next we turn to symmetries. A symmetry of a $3$-dimensional theory by an algebraic group $G$
induces an action of $G$ on its category of line operators (Section \ref{sec:gauging_lines}).
Given such an action on a category of line operators $A$ by
$\mathbb{E}_2$-monoidal functors, one can form its \emph{equivariantization} $A^G$, which
is an $\mathbb{E}_2$-algebra among $\Rep\left(G\right)$-modules, i.e.\ an object of
$\umor_2\left(\Mod_{\Rep\left(G\right)}\left(\Pre\right)\right)$.
This is the (de)equivariantization of \cite{DGNO} for finite groups and of \cite{AG:small}
for quantum groups, formulated using $1$-affineness of $BG$ \cite{G:ShvCat}.

\begin{introthm}[Proposition \ref{prop:equiv}, Remark \ref{rmk:inv_descent} and Lemma \ref{lem:geometric}]
\label{thm:intro_equiv}
Let $G$ be an affine algebraic group of finite type over a field of characteristic zero, acting on
an $\EE_2$-algebra $A$ in $\Pre$ by $\EE_2$-functors.
\begin{enumerate}
\item Equivariantization gives an equivalence
\begin{equation*}
\umor_2\left(\ShvCat\left(BG\right)\right)\simeq\umor_2\left(\Mod_{\Rep\left(G\right)}\left(\Pre\right)\right)
\ ,
\end{equation*}
and $A^G$ is invertible in the latter if and only if $A$ is invertible in $\umor_2\left(\Pre\right)$.
\item If the action comes from an action of $G$ on a perfect stack $X$ satisfying diagonal
descent, with $X/G$ perfect and
$A=\left(\QC\left(\cL X\right),*\right)$, then $\QC\left(X/G\right)$ trivializes
$A^G\simeq\QC\left(\cL X/G\right)$ over $\Rep\left(G\right)$.
\end{enumerate}
\end{introthm}

When $A$ is invertible, we define the \emph{anomaly} of the action as the class
$\left[A^G\right]-\left[\Rep\left(G\right)\otimes A\right]$ in the Witt group of
$\Rep\left(G\right)$-linear $\EE_2$-algebras (Definition \ref{defn:anomaly}). 
By (2), symmetries of the $2$-shifted cotangent bundle $T^*\left[2\right]X$ which are induced
by symmetries of $X$ itself are anomaly-free, at least when $X$ satisfies diagonal descent and
the line operators are $\QC\left(\cL X\right)$. This suggests that anomalies come from
symmetries of $T^*\left[2\right]X$ which do not come from $X$.
Two such symmetries motivate this paper. The first is the group of linear symmetries of
$T^*\left[2\right]\left(V\left[1\right]\right)\simeq\left(V\dsum V^\vee\right)\left[1\right]$, the
target of linear Rozansky--Witten theory, preserving its shifted symplectic structure; this is
larger than the group $\mathrm{GL}\left(V\right)$ of symmetries of $V\left[1\right]$.
The second is the quantum Frobenius symmetry of the small quantum group, to which we turn
next.

\subsection*{Quantum groups of tori}
\label{sec:intro_CS}

The line operators of Chern--Simons theory for a complex reductive algebraic group at a
root of unity are expected to form the $\EE_2$-category $\Rep\left(G_q\right)$ of
representations of the big quantum group.
Quantum Frobenius gives a central $\EE_2$-functor $\Rep\left(H\right)\to\Rep\left(G_q\right)$
for a dual group $H$, and $\Rep\left(G_q\right)$ is the equivariantization of an action of
$H$ on a version $A$ of the small quantum group \cite{AG:small, Neg:frob, Neg:log}; these
results are in the abelian setting (Remark \ref{rmk:abelian}).
For a torus, $A$ is a non-degenerate pointed braided fusion category, and for fusion categories
nothing is lost in passing to derived categories. So for tori we can work in the derived setting.

\begin{introthm}[Theorem \ref{thm:T_invertible}]
\label{thm:intro_C}
Let $k=\CC$, let $T$ be a torus with character lattice $\Lambda$, let $Q\colon\Lambda\to k^\times$ be a quadratic
form whose double braiding has radical $\Lambda_0$ of finite index with $Q|_{\Lambda_0}=1$, let
$\Rep\left(T_q\right)$ be the braided category of $\Lambda$-graded complexes determined by $Q$, and let
$H\ceqq\Hom\left(\Lambda_0,\GG_m\right)$. Then $\Rep\left(T_q\right)$ is the equivariantization of
an action of $H$ on the non-degenerate pointed braided fusion category $A=\Vect_{\Lambda/\Lambda_0}$,
and it is invertible in
\begin{equation*}
\umor_2\left(\Mod_{\Rep\left(H\right)}\left(\Pre\right)\right) \, .
\end{equation*}
\end{introthm}

Here $\Rep\left(T_q\right)$ and $A$ are the categories of representations of the big and small
quantum groups of $T$. The weight modules of Lusztig's divided-power quantum group of $T$ are
$\Lambda$-graded, with braiding given by $Q$, and $\Rep\left(H\right)=\Vect_{\Lambda_0}$ is the
M\"uger center of $\Rep\left(T_q\right)$. When $\Lambda_0$ is the character lattice of the target
of quantum Frobenius, its inclusion is quantum Frobenius, and $A=\Vect_{\Lambda/\Lambda_0}$ is the
category of representations of the small quantum group of $T$, which is the group algebra of the
finite group $\Hom\left(\Lambda/\Lambda_0,\GG_m\right)$ \cite{AG:small}. For instance, for
$T=\GG_m$ and $Q\left(a\right)=q^{a^2}$, with $q$ a primitive $n$-th root of unity and $n$ odd,
$\Lambda_0=n\ZZ$, quantum Frobenius is the isogeny $z\mapsto z^n$, and the small quantum group is
$k\left[K\right]/\left(K^n-1\right)$.

So the anomaly of the action of $H$ on $A$ is defined (Definition \ref{defn:anomaly}); we
do not determine here whether it vanishes. Detecting it, by transgression to character
stacks of surfaces, will be treated in forthcoming work.

\subsection*{Outlook: an analytic flavor}
In Section \ref{sec:aggregated} we briefly discuss an analytic setting, aggregated vector spaces, in
which more Witt classes should become trivial. We show that, given suitable classes attached to
real quadratic spaces, the image of a non-degenerate pointed braided fusion category there only depends on its
central charge modulo $8$ (Proposition \ref{prop:Witt_reduction}).
This analytic context, and various TQFTs defined in this setting, will be the subject of
forthcoming work with Pelle Steffens.

\subsection*{Organization} 

Section \ref{sec:Morita} recalls higher Morita categories and the criteria for
dualizability and invertibility. Section \ref{sec:sheaves} recalls quasi-coherent sheaves
on stacks and their factorization homology, and Section \ref{sec:inv_sheaves} studies
their invertibility. 
We introduce and prove our results concerning shifted $\cD$-modules in Section
\ref{sec:dmod}. Section \ref{sec:gauging} treats equivariantization and anomalies,
including quantum groups of tori and our expectations in an analytic setting.

\subsection*{Acknowledgments}

We thank David Ben-Zvi, Dan Freed, Pavel Safronov, Germ\'an Stefanich and Constantin Teleman
for many inspiring conversations. Special thanks go to Pelle Steffens for many useful
discussions, in particular on the analytic setting of Section \ref{sec:aggregated}.

Both authors were supported by the Simons Collaboration on Global Categorical Symmetries
(1013836 and 8528-03).
Both authors were supported by the Deutsche Forschungsgemeinschaft (DFG, German Research
Foundation) through the Collaborative Research Centers SFB 1085 Higher invariants -
224262486 and SFB 1624 Higher structures, moduli spaces and integrability - 506632645.

\section{Higher Morita categories and invertibility}
\label{sec:Morita}

\subsection{Conventions}
\label{sec:conventions}

Throughout, $k$ is a field of characteristic zero. We write $\Pre$ for the symmetric monoidal
$\left(\infty,2\right)$-category of presentable stable $k$-linear $\infty$-categories,
colimit-preserving $k$-linear functors and natural transformations, with the Lurie tensor product.
A stack is a derived stack over $k$; further hypotheses (perfect, QCA, bounded) are stated where
they are used. Formal completions and $\IndCoh_0$ are recalled in Section \ref{sec:IndCoh}.
We write $d$ for the dimension in $\umor_d$ and $n$ for the shift in $\DMod^n$.

Furthermore, $V$ always denotes a vector space and $V^\vee$ its dual. Shifts are written
$V\left[n\right]$ with the convention that $V\left[1\right]\simeq BV$; in particular it is an odd vector space, and
$T^*\left[2\right]\left(V\left[1\right]\right)\simeq\left(V\dsum V^\vee\right)\left[1\right]$.
Duals are always written
$\left(-\right)^\vee$; we reserve the notation $\left(-\right)^*$ for pullbacks, cotangent bundles and
(co)chains.

\begin{rmk}[Abelian versus derived]
\label{rmk:abelian}
Much of the literature on braided tensor categories that we use, e.g.\
\cite{BJS, BJSS, Neg:frob}, works in the $\left(2,2\right)$-category of locally
presentable $k$-linear (abelian) categories, whereas we work with presentable stable
$k$-linear $\infty$-categories (Section \ref{sec:conventions}). We pass from the former to the
latter by taking derived categories, and we only do this for fusion categories, where nothing is
lost: A finite semisimple category with simple objects $I$ becomes the product 
$\prod_{i\in I} \Vect$, and relative tensor products over fusion categories have no higher
derived terms, since module categories over a fusion category are semisimple. For
non-semisimple categories the derived relative tensor products and centers can differ from
the abelian ones, and we cite results about them only as statements in the abelian
setting. Results imported from the abelian setting are marked by a reference to this
remark.
\end{rmk}

\subsection{Morita category of \texorpdfstring{$\EE_d$}{Ed}-algebras}
\label{sec:mor}

Let $\coef$ be a symmetric monoidal $\infty$-cate\-gory which admits sifted colimits,
preserved by $\otimes$ in each variable, as in \cite{SSS}.
The Morita $\left(\infty,d+1\right)$-category $\umor_d\left(\coef\right)$ of \cite{SSS} has
$\EE_d$-algebras as objects, bimodules in $\EE_{d-1}$-algebras as $1$-morphisms, bimodules of
bimodules in $\EE_{d-2}$-algebras as $2$-morphisms, and so on, with homomorphisms as
$\left(d+1\right)$-morphisms. It is symmetric monoidal, and its loops are the Morita category
of $\EE_{d-1}$-algebras,
\begin{equation*}
\Omega\,\umor_d\left(\coef\right)\simeq\umor_{d-1}\left(\coef\right) \ ,
\end{equation*}
cf.\ \cite[Proposition 3.3.8]{Sch:thesis}, where this is shown for the
factorization-algebra model of $\umor_d$.
For $\coef$ a (symmetric monoidal) $\left(\infty,m\right)$-category, one
obtains a (symmetric monoidal) $\left(\infty,d+m\right)$-category \cite{JFS}. In particular, since $\Pre$ is a symmetric monoidal
$\left(\infty,2\right)$-category, $\umor_d\left(\Pre\right)$ is a symmetric monoidal
$\left(\infty,d+2\right)$-category.

In our examples, such as $\coef=\Pre$ and $\coef=\Mod_R\left(\Pre\right)$ for
$R\in\CAlg\left(\Pre\right)$, $\coef$ is not presentable, but it admits all small colimits,
preserved by $\otimes$ in each variable \cite[Section 4.8.1]{L:HA}.

\subsection{Dualizability of \texorpdfstring{$\mathbb{E}_d$}{Ed}-algebras}
\label{sec:dual_Ek}

By taking opposite algebra structures in each of the $d$ directions, every $\EE_d$-algebra is $d$-dualizable in $\umor_d\left(\coef\right)$, see \cite{GS}. For TQFT applications, the more interesting question is
 when it is $\left(d+1\right)$-dualizable, i.e.\ when it defines a fully extended
framed $\left(d+1\right)$-dimensional TQFT.
This is answered by the following theorem, conjectured by Lurie
\cite[Remark 4.1.27]{L:Cob}.

\begin{theorem}[{\cite[Theorem A]{SSS}}]
\label{thm:SSS_dual}
An $\EE_d$-algebra $A$ in $\coef$ is $\left(d+1\right)$-dualizable in $\umor_d\left(\coef\right)$ if,
and only if, for all $0\leq i\leq d$, the object $A$ is left dualizable as a left module over
the factorization homology
\begin{equation*}
\int_{S^{i-1}\times\RR^{d-i+1}}A \ .
\end{equation*}
\end{theorem}

For $i=0$ the factorization homology is $\int_{\emptyset}A = \one_{\coef}$, so the condition is
dualizability of $A$ in $\coef$.
Theorem \ref{thm:SSS_dual} also applies to Morita categories valued in a symmetric monoidal
$\left(\infty,m\right)$-category such as $\Pre$, since $\left(d+1\right)$-dualizability is
detected in the underlying $\left(\infty,d+1\right)$-category \cite[Section 1.5]{SSS}.

\begin{rmk}
\label{rmk:SSS_relative}
There is also a relative version for $1$-morphisms $M\colon A\to B$ in
$\umor_d\left(\coef\right)$, characterizing when $M$ is $d$-times right adjointable in terms of
dualizability over factorization homologies of $M$ \cite[Theorem B]{SSS}. Via \cite{JFS}, such
$M$ defines an oplax relative field theory. Two cases are relevant below. The regular module
$A$ always gives a maximally adjointable $1$-morphism. An $A$-central $\EE_{d-1}$-algebra $N$,
i.e.\ a $1$-morphism $A\to\one$
(Notation \ref{nota:E2}(4)), is
maximally adjointable if, and only if, $N$ is $d$-dualizable as an object of
$\umor_{d-1}\left(\coef\right)$ \cite[Section 4.5]{SSS}. For $d=2$, this is the dualizability input for the boundary
conditions of Section \ref{sec:line_ops}.
\end{rmk}

\subsection{Invertibility of \texorpdfstring{$\mathbb{E}_d$}{Ed}-algebras}
\label{sec:inv_Ek}

For an $\EE_d$-algebra $A$ and $0\leq j\leq d$, write
$\mathcal{Z}_j\left(A\right)\ceqq\End_{\int_{S^{j-1}\times\RR^{d-j+1}}A}\left(A\right)$ for its
$\EE_j$-center (cf.\ Notation \ref{nota:E2}). The action of
$\int_{S^{i-1}\times\RR^{d-i+1}}A$ on $A$ induces a canonical map to $\mathcal{Z}_{d-i}\left(A\right)$.

\begin{cor}[{\cite[Corollary A]{SSS}}]
\label{cor:SSS_inv}
An $\EE_d$-algebra $A$ in $\coef$ is invertible in $\umor_d\left(\coef\right)$ if, and only if, it
is $\left(d+1\right)$-dualizable and, for all $0\leq i\leq d$, the canonical map
\begin{equation*}
\int_{S^{i-1}\times\RR^{d-i+1}}A\longrightarrow \mathcal{Z}_{d-i}\left(A\right)
\end{equation*}
to the $\EE_{d-i}$-center is an equivalence.
\end{cor}

In particular, \cite{SSS} proves \cite[Conjecture 1.9]{BJSS}.
For $d=2$, the conditions are those of Definition \ref{defn:inv_E2} below:
non-degeneracy ($i=0$), factorizability ($i=1$) and cofactorizability ($i=2$).

\subsubsection{Invertibility of \texorpdfstring{$\mathbb{E}_2$}{E2}-algebras}
\label{sec:inv_E2}

This section unpacks Corollary \ref{cor:SSS_inv} for $d=2$.
We will make use of the usual description of $\mathbb{E}_2$-algebras as locally constant
factorization algebras on $\left(0,1\right)^2 \simeq \mathbb{R}^2$.

\begin{nota}
\label{nota:E2}
Let $A$ be an $\mathbb{E}_2$-algebra in $\coef$.
\begin{enumerate}
\item We write $A^{\rev}$ and $A^{\sigma\rev}$ for the pushforwards of $A$ along the
reflections of $\left(0,1\right)^2$ in the first and second coordinate. Rotating by $\pi$
provides equivalences $A^{\rev}\simeq A^{\sigma\rev}$ of $\mathbb{E}_2$-algebras, and
$A^{\rev}$ is the dual of $A$ in $\umor_2\left(\coef\right)$ \cite{GS}.
\item We write
\begin{equation*}
A^e\ceqq \int_{S^0\times\mathbb{R}\times\mathbb{R}}A\simeq A\otimes A^{\rev}\qquad\text{and}\qquad
\mathrm{HC}\left(A\right):=\int_{S^1\times\mathbb{R}}A \ ,
\end{equation*}
where $S^0\times\mathbb{R}\subset\mathbb{R}$ and $S^1\times\mathbb{R}\subset\mathbb{R}^2$
carry the framings bounding a disk, and $A\otimes A^{\sigma\rev}$ is obtained in the same
way from $\mathbb{R}\times S^0\times\mathbb{R}$. These are $\mathbb{E}_1$-algebras acting on $A$ via the radial coordinate. The notation $\mathrm{HC}\left(A\right)$ follows \cite{BJSS}, after the Harish-Chandra category of a braided tensor category.
\item The $\mathbb{E}_j$-centers of Section \ref{sec:inv_Ek} of an $\mathbb{E}_2$-algebra
$A$ are explicitly:
\begin{equation*}
\mathcal{Z}_0\left(A\right) \ceqq \End\left(A\right)  \ ,
\qquad
\mathcal{Z}_1\left(A\right) \ceqq \End_{A^e}\left(A\right)  \ , \quad\text{and}
\qquad
\mathcal{Z}_2\left(A\right) \ceqq \End_{\HC\left(A\right)}\left(A\right) \ .
\end{equation*}
\item A 1-morphism $A\rightarrow B$ in $\umor_2\left(\coef\right)$ is an $\mathbb{E}_1$-algebra in
$\left(A , B\right)$-bimodules, equivalently an $\mathbb{E}_1$-algebra $C$ in $\coef$
together with a map $A\otimes B^{\sigma\rev}\rightarrow \mathcal{Z}_1\left(C\right)$.
We call such data an $\left(A,B\right)$-central algebra. In particular, a $1$-morphism $B\to\one$ is a
$B$-central algebra, i.e.\ an $\EE_1$-algebra $C$ with a central structure $B\to \mathcal{Z}_1\left(C\right)$.
So, as in \cite{BJSS, SSS}, a boundary condition, or relative theory, for $B$ is a
$1$-morphism $B\to\one$. The reformulation of bimodules in terms of central algebras is
\cite[Section 5.3.1]{L:HA}.

\item For an $\mathbb{E}_1$-algebra $C$ we write $C^{!}\ceqq \mathrm{Hom}_{C^e}\left(C,C^e\right)$
and $C^\vee \ceqq \Hom\left(C , \one_{\coef}\right)$ for the duals of $C$
as a $C^e$-module and as an object of $\coef$, whenever these exist.
\end{enumerate}
\end{nota}

\begin{theorem}[{\cite[Theorem 2.26]{BJSS}}]
\label{thm:trivial_Alg2}
Let $\coef$ be a closed symmetric monoidal $\left(\infty,2\right)$-category admitting geometric
realizations, e.g.\ $\coef=\Pre$.
Let $B$ be an $\EE_2$-algebra in $\coef$ and let $C$ be a $B$-central algebra, viewed as a
1-morphism $C\colon B\to\one_{\coef}$ in $\umor_2\left(\coef\right)$. Then $C$ is invertible if,
and only if, $C$ is 2-dualizable in $\umor_1\left(\coef\right)$, that is, dualizable as an object of
$\coef$ and as a $C^e$-module, and the following maps are equivalences:
\begin{enumerate}
    \item the evaluation map $C^{!}\otimes_B C\to C^e$;
    \item the map $B\to \mathcal{Z}_1\left(C\right)$ given by the $B$-central structure on $C$;
    \item the evaluation map $C^{\vee}\otimes_{C\otimes_B C^{\rev}}C\to \one_{\coef}$;
    \item the map $C^{\rev}\otimes_B C\to \mathrm{End}\left(C\right)$ given by the left and
    right actions of $C$ on itself.
\end{enumerate}
\end{theorem}

Invertibility of $A$ itself is the case of Theorem \ref{thm:trivial_Alg2} in which
$B = A\otimes A^{\sigma\rev}$ and $C = A$, viewed as the evaluation $1$-morphism
$A\otimes A^{\sigma\rev}\to\one_{\coef}$ of the duality between $A$ and $A^{\sigma\rev}$, and the four conditions then assemble into the following three.

\begin{defn}[\cite{BJSS}]
\label{defn:inv_E2}
An $\mathbb{E}_2$-algebra $A$ in $\coef$ is
\begin{enumerate}
\item \emph{non-degenerate} if the unit map $\one_{\coef}\to \mathcal{Z}_2\left(A\right)$ is an
equivalence;
\item \emph{factorizable} if the action map $A\otimes A^{\sigma\rev}\to \mathcal{Z}_1\left(A\right)$ is an equivalence;
\item \emph{cofactorizable} if the action map
$\HC\left(A\right)\to \End\left(A\right)=\mathcal{Z}_0\left(A\right)$ is an equivalence.
\end{enumerate}
\end{defn}

\begin{theorem}[{\cite[Theorems 1.2 and 1.4]{BJSS};
\cite[Theorem A and Corollary A]{SSS}}]
\label{thm:inv_E2}
Let $A$ be an $\mathbb{E}_2$-algebra in $\coef$. Then $A$ is $3$-dualizable in
$\umor_2\left(\coef\right)$ if, and only if, $A$ is dualizable as an object of $\coef$, as a left
$A^e$-module, and as a left $\mathrm{HC}\left(A\right)$-module. Moreover, $A$ is invertible in
$\umor_2\left(\coef\right)$ if, and only if, it is $3$-dualizable, non-degenerate, factorizable and
cofactorizable.
\end{theorem}

This is the case $d=2$ of Theorem \ref{thm:SSS_dual} and Corollary \ref{cor:SSS_inv}; cf.\
\cite{BJS} for related results on braided tensor categories in the abelian setting.

\subsubsection{Braided tensor categories}
\label{sec:braided}

Braided tensor categories are $\mathbb{E}_2$-algebras in the $\left(2,2\right)$-category of
locally presentable $k$-linear categories. The results below are stated there; we only use them
in $\Pre$ for fusion categories (Remark \ref{rmk:abelian}).
In either setting, invertibility in the
Morita category $\umor_2$ is detected in its underlying $\left(\infty , 3\right)$-category
\cite[Section 1.5]{SSS}.

In Definition \ref{defn:inv_E2}, $\mathcal{Z}_1\left(A\right)$ is the Drinfeld center and, in the abelian
setting, $\mathcal{Z}_2\left(A\right)$ is the M\"uger center of $A$ \cite[Propositions 3.6 and 3.7]{BJSS}, so non-degeneracy is the
familiar condition that the M\"uger center is trivial. Over an algebraically closed field
of characteristic zero, the three conditions of Definition \ref{defn:inv_E2} are mutually equivalent
for finite braided tensor categories, and a finite braided tensor category is invertible
if and only if it is non-degenerate \cite[Theorems 1.1 and 3.20]{BJSS}.

Examples include non-semisimple categories such as representation categories of small
quantum groups at suitable roots of unity, which give invertible fully extended framed
$4$-dimensional topological field theories in the abelian setting. Outside the finite setting this is not true:
The category of integrable representations of a quantum group at a generic
parameter is non-degenerate but not factorizable \cite[Remark 1.7]{BJSS}.

\subsection{Invertibility of commutative algebras}
\label{sec:inv_comm}

For a commutative algebra, invertibility in the higher Morita category is equivalent to
triviality.

\begin{prop}
\label{prop:inv_symm}
Let $A\in\CAlg\left(\coef\right)$ and $d\geq 1$.
Then $A$ is invertible in $\umor_d\left(\coef\right)$ if and only if the unit
$u\colon\one_{\coef}\to A$ is an equivalence.
\end{prop}

\begin{proof}
The ``if'' direction is clear, since then $A\simeq\one_{\coef}$ is the monoidal unit.
For ``only if'' we use the following observation.
Suppose we are given maps
$\one_{\coef}\xrightarrow{u}A\xrightarrow{i}W\xrightarrow{r}A$
with $ri\simeq\Id_A$ and $iu$ an equivalence.
Then $i$ is split epi, via $u\left(iu\right)^{-1}$, and split mono, via $r$, so $i$ is an
equivalence.
Hence $u\simeq i^{-1}\left(iu\right)$ is an equivalence.

Write $R_d\ceqq\int_{S^{d-1}\times\RR}A$, so that
$\mathcal{Z}_d\left(A\right)=\End_{R_d}\left(A\right)$.
If $A$ is invertible, then by Corollary \ref{cor:SSS_inv}, applied to $\int_{\emptyset}A=\one_{\coef}$, the unit
$\one_{\coef}\to\mathcal{Z}_d\left(A\right)$, which picks out $\Id_A$, is an equivalence.
Since $A$ is commutative, $R_d$ is the tensor $S^{d-1}\otimes A$ in $\CAlg\left(\coef\right)$
\cite[Section 5.5]{L:HA}, and it acts on $A=\int_{\RR^d}A$ through the map of commutative
algebras $\varepsilon\colon R_d\to A$ induced by $S^{d-1}\to\pt$.
Let $i\colon A\to\End_{R_d}\left(A\right)$ send $a$ to multiplication by $a$; this is $R_d$-linear
since $A$ is commutative and $R_d$ acts through $\varepsilon$.
By the unit axiom, $i\left(1_A\right)\simeq\Id_A$, so the unit of $\End_{R_d}\left(A\right)$ is $iu$.
Evaluation at $1_A$ defines $r\colon\End_{R_d}\left(A\right)\to A$ with $ri\simeq\Id_A$, and we
conclude with the above observation.
For $d=1$, $R_1=A\otimes A^{\rev}$ and this condition is centrality in \cite[Theorem 1.3]{BJSS}.
\end{proof}

Applying Proposition \ref{prop:inv_symm} with $\coef$ replaced by
$\Mod_R\left(\coef\right)$, whose unit is $R$, gives the relative version.

\begin{cor}
\label{cor:inv_symm_rel}
Let $R\in\CAlg\left(\coef\right)$, $A\in\CAlg\left(\Mod_R\left(\coef\right)\right)$ and
$d\geq 1$.
Then $A$ is invertible in $\umor_d\left(\Mod_R\left(\coef\right)\right)$ if and only if the
unit $R\to A$ is an equivalence.
\end{cor}

\subsection{Line operators in \texorpdfstring{$3$D}{three-dimensional} TQFTs}
\label{sec:line_ops}

Let $\cT$ be a symmetric monoidal $\left(\infty,3\right)$-category and
$\mathcal{F}\colon\Bord_3^{\mathrm{fr}}\to\cT$ a framed fully extended TQFT, with $B\ceqq
\mathcal{F}\left(\pt\right)$.
The \emph{category of line operators} of $\mathcal{F}$ is $\mathcal{F}\left(S^1\right)$, where $S^1$
carries the framing it inherits as the boundary of a disk; it is an $\EE_2$-algebra via embeddings
of disks.
By the cobordism hypothesis \cite{L:Cob},
\begin{equation*}
\mathcal{F}\left(S^1\right)\simeq\End_{\End_\cT\left(B\right)}\left(\Id_B\right) \ .
\end{equation*}
For $\cT=\umor_1\left(\Pre\right)$, $B$ is an $\EE_1$-algebra in $\Pre$,
$\End_\cT\left(B\right)$ is the category of $B$-bimodules and $\Id_B = B$. So the line
operators form the Drinfeld center
\begin{equation*}
\mathcal{F}\left(S^1\right)\simeq \mathcal{Z}_1\left(B\right)=\End_{B^e}\left(B\right)
\end{equation*}
of Notation \ref{nota:E2}.
Moreover $B$ is tautologically $\mathcal{Z}_1\left(B\right)$-central, i.e.\ it is a $1$-morphism
$\mathcal{Z}_1\left(B\right)\to\one$ in $\umor_2\left(\Pre\right)$, which is maximally adjointable since
$B$ is fully dualizable in $\umor_1\left(\Pre\right)$ (Remark \ref{rmk:SSS_relative}). So $\mathcal{F}$ is a
boundary condition for the $4$-dimensional theory determined by
$\mathcal{Z}_1\left(B\right)$ (when that theory exists).

\section{Sheaves on stacks}
\label{sec:sheaves}

\subsection{Quasi-coherent sheaves}
\label{sec:QC}

To any stack $X$, there is an associated stable $\infty$-category $\QC\left(X\right)$ of
quasi-coherent sheaves.
The homotopy category of this is the usual unbounded quasi-coherent derived category.
By definition, $\QC$ is defined on affines $\Spec\left(A\right)$ to be $\Mod_A$.
The definition is extended to a general stack via right Kan extension:
Every stack $X$ can be written as a colimit of affines, and $\QC\left(X\right)$ is then
defined to be the limit of the corresponding diagram in $\Pre$:
\begin{equation*}
\QC\left(X\right)\ceqq \lim_{\Spec A \to X} \Mod_A \in \Pre
\end{equation*}

\subsubsection{Perfect stacks}

We will occasionally make the mild assumption that our stack is well-behaved in the
following sense.

\begin{defn}[{\cite[\S 3.1]{BZFN}}]
\label{defn:perfect}
A stack $X$ is \emph{perfect}
if it has affine diagonal and
$\QC\left(X\right)\simeq \Ind\Perf\left(X\right)$.
\end{defn}

This is a very broad class of stacks, including quasi-compact derived schemes with affine
diagonal, and linear quotients of quasi-projective derived schemes in characteristic zero
\cite[\S 3.3]{BZFN}. In particular, many instances of classifying stacks are perfect.

\begin{theorem}[{\cite[\S 4.2-4.3]{BZFN}}]
\label{thm:BZFN}
Consider maps $X\xrightarrow{p} Y \xleftarrow{p'} X'$ of perfect stacks.
Then $\QC\left(X\right)$ (and $\QC\left(X'\right)$) are self-dual over
$\QC\left(Y\right)$, and there are equivalences
\begin{equation*}
\QC\left(X\right) \otimes_{\QC\left(Y\right)} \QC\left(X'\right) \simeq
\QC\left(X\times _Y X'\right) \simeq
\Fun_{\QC\left(Y\right)}\left(\QC\left(X\right) , \QC\left(X'\right)\right)
\end{equation*}
of $\QC\left(X\right)$-$\QC\left(X'\right)$-bimodule categories.
\end{theorem}

\begin{rmk}
\label{rmk:BZFN_monoidal}
The former equivalence is symmetric monoidal, i.e. an equivalence in
$\CAlg\left(\Pre\right)$.

In the case $X = X'$, the functors from $\QC\left(X\right)$ to itself naturally form a
monoidal category with respect to composition. When considered with this monoidal
structure, the latter equivalence also turns out
to be monoidal, but where the monoidal structure on $\QC\left(X\times_Y X\right)$ is taken
to be convolution (see Section \ref{sec:convolution}).
Explicitly, the latter equivalence sends a kernel $K\in \QC\left(X\times_Y X'\right)$ to
the integral transform $\pi'_{*} \left( \pi^* \left(-\right)\otimes K\right)$, where
$\pi\colon X\times_Y X'\to X$ and $\pi'\colon X\times_Y X'\to X'$ are the projections.
\end{rmk}

\subsubsection{Symmetric-monoidal structure}

For a fixed perfect stack $S$, let $\Stk_S^{\perf}$ be the $\infty$-category of perfect
stacks over $S$, and write
$\Mod_{\QC\left(S\right)}\ceqq\Mod_{\QC\left(S\right)}\left(\Pre\right)$.
Pullback along $X\to S$ makes $\QC\left(X\right)$ a commutative algebra in
$\Mod_{\QC\left(S\right)}$, so we obtain a functor
\begin{equation}
\QC\colon\left(\Stk_S^{\perf}\right)^{\op}\longrightarrow\CAlg\left(\Mod_{\QC\left(S\right)}\right)
\ .
\label{eqn:QC_CAlg}
\end{equation}
Pushouts in $\CAlg\left(\Mod_{\QC\left(S\right)}\right)$ are relative tensor products
\cite[Proposition 3.2.4.7 and Section 4.5]{L:HA}, and $\Stk_S^{\perf}$ is closed under fiber products
\cite[\S 3]{BZFN}.
So Theorem \ref{thm:BZFN} says precisely that \eqref{eqn:QC_CAlg} sends finite limits to
finite colimits.

\subsection{Factorization homology}
\label{sec:fac_hom}

\subsubsection{Viewing the commutative object in the higher Morita categories}
\label{sec:QC_Algn}

Pointwise tensor product of sheaves endows $\QC\left(X\right)$ with a symmetric monoidal structure, so
a stack $X$ defines an object
$\left(\QC\left(X\right) , \otimes\right)\in\umor_d\left(\Pre\right)$ for all $d$.
Its factorization homology is computed by mapping stacks.

\begin{theorem}[{\cite[\S 4]{BZFN}; see \cite{Ste:Tan} for weaker hypotheses}]
\label{thm:int_QC}
Let $X$ be a perfect stack, and let $M$ be a smooth $m$-manifold, $m\leq d$, which is the
interior of a compact manifold with boundary (so that $M$ has finite homotopy type). Since
$\left(\QC\left(X\right),\otimes\right)$ is commutative, its factorization homology over $M$ is
defined without a framing, and agrees with factorization homology of the $\EE_d$-algebra
$\QC\left(X\right)$ over $M\times\RR^{d-m}$. There is a natural symmetric monoidal
equivalence
\begin{equation}
\int_M \QC\left(X\right) \simeq \QC\left(X^M\right) \ ,
\qquad X^M\ceqq\Maps\left(M_B , X\right) \ ,
\label{eqn:int_QC}
\end{equation}
where $M_B$ is the constant (Betti) stack on $M$ and $\Maps$ is the derived mapping stack.
\end{theorem}

\begin{proof}
Since $\left(\QC\left(X\right),\otimes\right)$ is commutative, $\int_M\QC\left(X\right)$ only
depends on the homotopy type of $M$: It is the tensor $M\otimes\QC\left(X\right)$ in
$\CAlg\left(\Pre\right)$, i.e.\ the colimit of the constant diagram with value
$\QC\left(X\right)$ indexed by $M$ \cite[Section 5.5]{L:HA}.
As $M$ is a finite space, this colimit is built from copies of $\QC\left(X\right)$ by finitely
many pushouts, which in $\CAlg\left(\Pre\right)$ are relative tensor products.
Dually, $X^M$ is the corresponding finite limit of copies of $X$, which is again perfect
\cite[\S 3]{BZFN}.
By Theorem \ref{thm:BZFN}, the functor \eqref{eqn:QC_CAlg} sends these finite limits to the
finite colimits computing $M\otimes\QC\left(X\right)$.
\end{proof}

Applying Theorem \ref{thm:BZFN} with $Y = X^2$ (equivalently, Theorem \ref{thm:int_QC} with
$M = S^1$), we obtain the following.

\begin{theorem}[\cite{BZFN}]
\label{thm:BZFN_loops}
For a perfect stack $X$, there are equivalences
\begin{equation*}
\QC\left(X\right) \otimes_{\QC\left(X^2\right)} \QC\left(X\right)
\simeq \QC\left(\mathcal{L}X\right) \simeq
\mathcal{Z}_1\left(\QC\left(X\right) \right) \ .
\end{equation*}
\end{theorem}

\subsection{Dualizability}
\label{sec:QC_dualizability}

\begin{prop}
\label{prop:QC_dualizable}
Let $X$ be a perfect stack. Then $\left(\QC\left(X\right),\otimes\right)$ is
$\left(d+1\right)$-dualizable in $\umor_d\left(\Pre\right)$ for every $d\geq0$.
\end{prop}

\begin{proof}
By Theorem \ref{thm:SSS_dual}, we need $\QC\left(X\right)$ to be dualizable as a module over
\begin{equation*}
\int_{S^{i-1}\times\RR^{d-i+1}}\QC\left(X\right)
\end{equation*}
for $0\leq i\leq d$. Since
$\QC\left(X\right)$ is commutative, this is $\QC\left(X^{S^{i-1}}\right)$ by Theorem
\ref{thm:int_QC}, acting via pullback along the constant maps $X\to X^{S^{i-1}}$ (cf.\ the
proof of Proposition \ref{prop:inv_symm}). The stack $X^{S^{i-1}}$ is a finite limit of copies
of $X$, hence perfect \cite[\S 3]{BZFN}, so $\QC\left(X\right)$ is self-dual over it by Theorem
\ref{thm:BZFN}.
\end{proof}

The argument uses that $\QC\left(X\right)$ is commutative. It does not apply to convolution
categories such as $\left(\QC\left(\cL X\right),*\right)$ or $\DMod^2\left(X\right)$; for the
latter see Corollary \ref{cor:3dual}.

\subsection{The convolution \texorpdfstring{$\EE_2$}{E2}-structure}
\label{sec:LX_convolution}

For a perfect stack $X$, we equip $\QC\left(\cL X\right)$ with the $\EE_2$-monoidal structure
of the Drinfeld center $\mathcal{Z}_1\left(\QC\left(X\right)\right)$, via Theorem
\ref{thm:BZFN_loops}, and write $\left(\QC\left(\cL X\right),*\right)$ for this
$\EE_2$-algebra. Its underlying monoidal structure is convolution relative to $X^2$ (Section
\ref{sec:convolution} with $Y=X^2$): Composition in $\mathcal{Z}_1$ is convolution on $\cL X$.

\section{Invertibility of categories of sheaves}
\label{sec:inv_sheaves}

We now ask when categories of sheaves are invertible in the higher Morita categories of Section
\ref{sec:Morita}. For the symmetric monoidal structure the answer is essentially never (Section
\ref{sec:inv_QC}). We therefore turn to convolution structures (Section \ref{sec:convolution}),
and in particular to quasi-coherent sheaves on loop stacks, which are Morita trivial whenever the
stack satisfies diagonal descent (Section \ref{sec:diag_descent}). This includes the formal
completion of a vector space at the origin (Section \ref{sec:formal_vector_space}).

\subsection{Symmetric monoidal categories of sheaves}
\label{sec:inv_QC}

We work relative to a base stack $S$; the absolute case is $S=\pt$. For a map $f\colon X\to S$,
the unit of the commutative $\QC\left(S\right)$-algebra $\left(\QC\left(X\right),\otimes\right)$ is
$f^*\colon\QC\left(S\right)\to\QC\left(X\right)$, so by Corollary \ref{cor:inv_symm_rel} it is
invertible in $\umor_d\left(\Mod_{\QC\left(S\right)}\left(\Pre\right)\right)$ if and only if $f^*$ is
an equivalence. To translate this into a condition on stacks, we use Tannaka duality.
\begin{lem}
\label{lem:Tannaka}
Let $f\colon X\to Y$ be a map of stacks such that $f^*$ is an equivalence.
\begin{enumerate}
\item If $f$ has a retraction $g$, i.e.\ $gf\simeq\Id_X$, and $Y$ is quasi-compact with
affine diagonal, then $f$ is an equivalence.
\item If $f$ has a section $s$, i.e.\ $fs\simeq\Id_Y$, and $X$ is quasi-compact with
affine diagonal, then $f$ is an equivalence.
\item If $f$ is flat and $X$ and $Y$ are quasi-compact with affine diagonal, then $f$ is an
equivalence.
\end{enumerate}
\end{lem}

\begin{proof}
For (1): $f^*g^*\simeq\Id$, so $g^*\simeq\left(f^*\right)^{-1}$, and hence
$\left(fg\right)^*\simeq g^*f^*\simeq\Id$ as symmetric monoidal functors.
By Tannaka duality \cite[Theorem 9.2.0.2]{L:SAG}, \cite{BHL}, maps into $Y$
are detected by their pullback functors, so $fg\simeq\Id_Y$.
Part (2) is the same argument with the roles reversed: $s^*\simeq\left(f^*\right)^{-1}$,
so $\left(sf\right)^*\simeq\Id$ and $sf\simeq\Id_X$.

For (3): Since $f$ is flat, $f^*$ is t-exact, and being conservative it detects connective
and discrete objects: If $f^*E$ is connective, then
$f^*\tau_{<0}E\simeq\tau_{<0}f^*E\simeq 0$, so $\tau_{<0}E\simeq 0$, and similarly for
$\tau_{>0}$.
Let $F\ceqq\left(f^*\right)^{-1}$. Then $F$ preserves connective objects, since
$f^*F\simeq\Id$. It also preserves flat objects, using that a connective object $E$ is flat if and only if
$E\otimes M$ is discrete for every discrete $M$ (affine-locally, \cite[Section 7.2.2]{L:HA}): If $E$ is flat
and $M\in\QC\left(Y\right)$ is discrete, then $f^*\left(F\left(E\right)\otimes M\right)\simeq
E\otimes f^*M$ is discrete, hence so is $F\left(E\right)\otimes M$.
So $F$ is a colimit-preserving symmetric monoidal functor preserving connective and flat
objects, and by Tannaka duality it is $h^*$ for some $h\colon Y\to X$.
Then $\left(fh\right)^*\simeq h^*f^*\simeq\Id$ and $\left(hf\right)^*\simeq f^*h^*\simeq\Id$,
so $fh\simeq\Id_Y$ and $hf\simeq\Id_X$.
\end{proof}

\begin{theorem}
\label{thm:inv_rel}
Let $f\colon X\to S$ be a map of quasi-compact stacks with affine diagonal, and let
$d\geq 1$.
Suppose that $f$ is flat, or that $f$ has a section, or that $f$ has a retraction.
Then the commutative $\QC\left(S\right)$-algebra $\left(\QC\left(X\right),\otimes\right)$
is invertible in $\umor_d\left(\Mod_{\QC\left(S\right)}\left(\Pre\right)\right)$ if and only if
$f$ is an equivalence.
\end{theorem}

\begin{proof}
By Corollary \ref{cor:inv_symm_rel}, invertibility is equivalent to $f^*$ being an
equivalence, and by Lemma \ref{lem:Tannaka} this is equivalent to $f$ being an equivalence.
\end{proof}

\subsection{Convolution categories and descent}
\label{sec:convolution}
\label{sec:conv}

When the stack on which we study sheaves is a fiber product, convolution gives
noncommutative monoidal structures.

Given a map of stacks $X\to Y$, the category of sheaves on the fiber product $X\times_Y X$
has a monoidal structure given by convolution.
Namely, we have a diagram of stacks:
\begin{equation*}
\begin{tikzcd}[column sep=small]
 & X\times_Y X\times_Y X
   \arrow[dl, "p_{12}"']
   \arrow[d, "p_{13}"]
   \arrow[dr, "p_{23}"] & \\
 X\times_Y X & X\times_Y X & X\times_Y X
\end{tikzcd}
\end{equation*}
and convolution
is induced by pulling and pushing along this diagram:
\begin{equation*}
M* N = p_{13*}\left(p_{12}^*M\otimes p_{23}^*N\right)
\end{equation*}

\begin{exm}
When $Y = \pt$, $\left(\QC\left(X\times X\right),*\right)$ is the monoidal category of integral
kernels, equivalent for $X$ perfect to $\End\left(\QC\left(X\right)\right)$ under composition
(Remark \ref{rmk:BZFN_monoidal}).
\end{exm}

\begin{theorem}[{\cite[Theorem 1.3]{BZFN:app}}]
Let $\pi \colon X \to Y$ be a map of perfect stacks.
If either
\begin{itemize}
\item $\pi$ is faithfully flat, or
\item $X$ and $Y$ are smooth and $\pi$ is proper and surjective,
\end{itemize}
then tensoring with
the $\QC\left(Y\right)$-$\QC\left(X\times_Y X\right)$-bimodule $\QC\left(X\right)$ defines
an equivalence
\begin{equation*}
\Mod_{\left(\QC\left(Y\right) , \otimes\right)}\left(\Pre\right) \lto{\sim}
\Mod_{\left(\QC\left(X\times_Y X\right) , *\right)}\left(\Pre\right)
\end{equation*}
\label{thm:BZFN_app}
\end{theorem}

\subsection{Quasi-coherent sheaves on loop stacks}

A particularly important instance of a convolutional category of sheaves on a fiber
product is given by the diagonal map $X \to Y = X\times X$.
This is the category obtained via Theorem \ref{thm:BZFN_loops} as the Drinfeld center of
the symmetric-monoidal category $\QC\left(X\right)$.

\subsubsection{Descent along the diagonal}
\label{sec:diag_descent}

\begin{defn}
\label{defn:diag_descent}
A stack $X$ satisfies \emph{diagonal descent} if
\begin{enumerate}
\item the conclusions of Theorem \ref{thm:BZFN} hold for the diagonal: $\QC\left(X\right)$ is
self-dual in $\Pre$ and over $\QC\left(X^2\right)$,
$\QC\left(X\right)\otimes\QC\left(X\right)\simeq\QC\left(X^2\right)$, and
\begin{equation*}
\QC\left(X\right)\otimes_{\QC\left(X^2\right)}\QC\left(X\right)\simeq\QC\left(\cL X\right)
\simeq\End_{\QC\left(X^2\right)}\left(\QC\left(X\right)\right) \ ,
\end{equation*}
where the second equivalence identifies composition with convolution; and
\item the conclusion of Theorem \ref{thm:BZFN_app} holds for the diagonal, i.e.\ tensoring
with $\QC\left(X\right)$ gives an equivalence
\begin{equation*}
\Mod_{\left(\QC\left(X^2\right),\otimes\right)}\left(\Pre\right)\simeq
\Mod_{\left(\QC\left(\cL X\right),*\right)}\left(\Pre\right) \ .
\end{equation*}
\end{enumerate}
\end{defn}

Condition (1) holds for every perfect stack by Theorem \ref{thm:BZFN} and Remark
\ref{rmk:BZFN_monoidal}, and it implies $\mathcal{Z}_1\left(\QC\left(X\right)\right)\simeq\QC\left(\cL X\right)$
as in Theorem \ref{thm:BZFN_loops}.
By Theorem \ref{thm:BZFN_app}, a perfect stack satisfies (2), and hence diagonal descent,
if $\Delta\colon X\to X\times X$ is faithfully flat, or if $X$ is smooth and $\Delta$ is
proper and surjective.

\begin{exm}
\label{exm:BG_diag_descent}
If $G$ is an affine algebraic group of finite type over a field of characteristic zero, then
$BG$ is perfect \cite[\S 3.3]{BZFN} and satisfies diagonal descent.
Indeed, the base change of $\Delta$ along the cover $\pt\to BG\times BG$ is
$\left(G\times G\right)/G\simeq G\to\pt$, so $\Delta$ is faithfully flat.
\end{exm}

\begin{rmk}
\label{rmk:diag_descent_fails}
Diagonal descent fails for every quasi-compact separated scheme $X$ with more than one point,
e.g.\ any positive-dimensional scheme of finite type over $k$.
Indeed, $X$ is perfect, and its diagonal is not surjective, so there is a non-empty affine open
$U\subset X^2$ disjoint from $\Delta\left(X\right)$. By Theorem \ref{thm:BZFN},
\begin{equation*}
\QC\left(X\right)\otimes_{\QC\left(X^2\right)}\QC\left(U\right)\simeq
\QC\left(X\times_{X^2}U\right)\simeq\QC\left(\emptyset\right)=0 \ ,
\end{equation*}
whereas $\QC\left(U\right)\neq0$. So tensoring with $\QC\left(X\right)$ is not conservative,
and condition (2) of Definition \ref{defn:diag_descent} fails.
\end{rmk}

The formal completion of a vector space at the origin also satisfies diagonal descent,
although it is not perfect and its diagonal is not flat (Proposition
\ref{prop:Morita_formal}).
Stacks satisfying diagonal descent are exactly those to which the Morita triviality result,
Theorem \ref{thm:trivial_QC}, applies.

\subsubsection{Invertibility}

\begin{lem}[Morita context]
\label{lem:Morita_context}
Let $X$ be a stack satisfying condition (1) of Definition \ref{defn:diag_descent}, and write
$R\ceqq\QC\left(X^2\right)$, $M\ceqq\QC\left(X\right)$ and
$A\ceqq\left(\QC\left(\cL X\right),*\right)\simeq\End_R\left(M\right)$. Then
$\Psi\ceqq M\otimes_A-\colon\Mod_A\left(\Pre\right)\to\Mod_R\left(\Pre\right)$ is left adjoint
to $\Phi\ceqq M^\vee\otimes_R-$ and fully faithful, and
\begin{equation*}
P\ceqq M\otimes_A M^\vee=\QC\left(X\right)\otimes_{\QC\left(\cL X\right)}\QC\left(X\right)
\lto{\ \varepsilon\ } R
\end{equation*}
is an idempotent $R$-coalgebra, i.e.\ $P\otimes_R P\simeq P$ via $\varepsilon$. The functor $\Phi$ is an equivalence if and only if the counit
$\varepsilon$ is an equivalence, which is condition (2) of Definition \ref{defn:diag_descent}.
\end{lem}

\begin{proof}
Since $M$ is dualizable over $R$ with dual $M^\vee\simeq M$, the right adjoint of $\Psi$ is
$\Fun_R\left(M,-\right)\simeq\Phi$, and $M^\vee\otimes_R M\simeq A$ as $A$-bimodules, so the
unit $\Id\to\Phi\Psi$ is an equivalence. In the other order, $\Psi\Phi\simeq P\otimes_R-$, with
counit $\varepsilon\otimes_R-$, and
$P\otimes_R P\simeq M\otimes_A\left(M^\vee\otimes_R M\right)\otimes_A M^\vee\simeq P$.
\end{proof}

The following upgrades the Morita equivalence of \cite{BZFN:app} in $\umor_1\left(\Pre\right)$
(Theorem \ref{thm:BZFN_app}) to $\umor_2\left(\Pre\right)$, under the different hypothesis of diagonal
descent, using \cite{BZFN} and the criterion of \cite[Theorem 2.26]{BJSS}.

\begin{theorem}
\label{thm:trivial_QC}
Let $X$ be a stack satisfying diagonal descent (Definition \ref{defn:diag_descent}), e.g.\
$X = BG$ (Example \ref{exm:BG_diag_descent}) or the formal completion $\widehat V_0$ of a
vector space at the origin (Proposition \ref{prop:Morita_formal}).
Then $\left(\QC\left(\mathcal{L}X\right),*\right)$ is Morita trivial, i.e.\ equivalent to the
unit of $\umor_2\left(\Pre\right)$.
More precisely, $\QC\left(X\right)$, regarded as a $\QC\left(\cL X\right)$-central algebra via
$\mathcal{Z}_1\left(\QC\left(X\right)\right)\simeq\QC\left(\cL X\right)$, is an invertible $1$-morphism
from $\QC\left(\mathcal{L}X\right)$ to the monoidal unit of $\umor_2\left(\Pre\right)$.
\end{theorem}

\begin{proof}
We check the criterion of Theorem \ref{thm:trivial_Alg2} \cite[Theorem 2.26]{BJSS}, for
$B = \left(\QC\left(\cL X\right),*\right)$ and $C = \left(\QC\left(X\right),\otimes\right)$. Note
that $C^e = C\otimes C^\rev\simeq\left(\QC\left(X^2\right),\otimes\right)$ since $C$ is
commutative.
By condition (1) of Definition \ref{defn:diag_descent}, $C$ is dualizable in $\Pre$ and over
$C^e$, i.e.\ $2$-dualizable in $\umor_1\left(\Pre\right)$, and $C$ is self-dual in both senses,
so $C^!\simeq C^\vee\simeq\QC\left(X\right)$.
\begin{enumerate}
\item The evaluation $C^!\otimes_B C\to C^e$ is the counit
$\QC\left(X\right)\otimes_{\QC\left(\cL X\right)}\QC\left(X\right)\to\QC\left(X^2\right)$ of
Lemma \ref{lem:Morita_context}. It is an equivalence
by condition (2).
\item The map $B\to \mathcal{Z}_1\left(C\right)$ is the equivalence of condition (1), cf.\ Theorem
\ref{thm:BZFN_loops}.
\item[(4)] Since $C$ is commutative, $C^\rev\otimes_B C\to\End\left(C\right)$ is the same counit,
followed by the equivalence $\QC\left(X^2\right)\simeq\Fun\left(\QC\left(X\right),\QC\left(X\right)\right)$
given by self-duality of $\QC\left(X\right)$ in $\Pre$. So it is an equivalence.
\item[(3)] By (4), $C\otimes_B C^\rev\simeq\End\left(C\right)$, and the evaluation in (3)
becomes $C^\vee\otimes_{\End\left(C\right)}C\to\Vect$. This is an equivalence for any dualizable
$C\in\Pre$, since $\End\left(C\right)\simeq C\otimes C^\vee$ is Morita trivialized by $C$. \qedhere
\end{enumerate}
\end{proof}

In particular, by Theorem \ref{thm:inv_E2}, $\left(\QC\left(\cL X\right),*\right)$ is
non-degenerate, factorizable and cofactorizable. For instance, for $X=BG$, where $\cL BG\simeq G/G$ is the adjoint quotient, we have
\begin{equation*}
\mathcal{Z}_2\left(\QC\left(G/G\right),*\right)\simeq\Vect
\qquad\text{and}\qquad
\QC\left(G/G\times G/G\right)\simeq \mathcal{Z}_1\left(\QC\left(G/G\right),*\right) \ ,
\end{equation*}
the latter via the action map, generalizing the modularity of the Drinfeld center of
$\Rep\left(G\right)$ for finite $G$.

\subsubsection{The formal completion of a vector space at the origin}
\label{sec:formal_vector_space}

Let $V$ be a finite-dimensional vector space and $\widehat V_0\ceqq V^\wedge_0$ its formal
completion at the origin. Strictly speaking $\widehat V_0$ is an ind-scheme rather than a
perfect stack, so Theorems \ref{thm:BZFN} and \ref{thm:BZFN_app} do not literally apply to
it. Nonetheless it satisfies diagonal descent.

\begin{prop}
\label{prop:Morita_formal}
The formal completion $\widehat V_0$ satisfies diagonal descent
(Definition \ref{defn:diag_descent}). In particular, tensoring with $\QC\left(\widehat V_0\right)$
gives a Morita equivalence between $\left(\QC\left(\widehat V_0^2\right),\otimes\right)$ and
$\left(\QC\left(\cL\widehat V_0\right),*\right)$, and
$\mathcal{Z}_1\left(\QC\left(\widehat V_0\right)\right)\simeq\QC\left(\cL\widehat V_0\right)$.
\end{prop}

\begin{proof}
\emph{Step 0: Completions as supports.}
For a derived scheme $Y$ of finite type and a closed point $0\in Y$, local duality identifies
$\QC\left(Y^\wedge_0\right)$, symmetric monoidally, with the category $\QC_0\left(Y\right)$ of
quasi-coherent sheaves set-theoretically supported at $0$
\cite[Chapters 7--8]{L:SAG}. The local cohomology functor $\Gamma_0$ is
smashing, $\Gamma_0 F\simeq F\otimes\Gamma_0\cO_Y$, so $\QC_0\left(Y\right)$ is an idempotent
commutative $\QC\left(Y\right)$-algebra, and for $f\colon Z\to Y$ we have
$\QC\left(Z\right)\otimes_{\QC\left(Y\right)}\QC_0\left(Y\right)\simeq\QC_{f^{-1}\left(0\right)}\left(Z\right)$.
Completion commutes with fiber products, so $\widehat V_0^n\simeq\left(V^n\right)^\wedge_0$ and
$\cL\widehat V_0\simeq\left(\cL V\right)^\wedge_0=\widehat V_0\times V\left[-1\right]$.
Hence $\QC\left(\widehat V_0^n\right)\simeq\QC_0\left(V^n\right)$ and
$\QC\left(\cL\widehat V_0\right)\simeq\QC_0\left(\cL V\right)$.

\emph{Step 1: Condition (1).}
Write $E\ceqq\QC_0\left(V^2\right)$, $M\ceqq\QC_0\left(V\right)$ and
$A\ceqq\QC_0\left(\cL V\right)$. Since the origin of $V^2$ has preimage the origin under both
$\Delta\colon V\to V^2$ and $\cL V\to V^2$, Step 0 gives
$M\simeq\QC\left(V\right)\otimes_{\QC\left(V^2\right)}E$ and
$A\simeq\QC\left(\cL V\right)\otimes_{\QC\left(V^2\right)}E$; moreover $M\otimes M\simeq E$.
Base change along the symmetric monoidal functor $\QC\left(V^2\right)\to E$ preserves duals,
relative tensor products, and endomorphisms of dualizable modules. So Theorem
\ref{thm:BZFN} for $\Delta\colon V\to V^2$, together with Remark \ref{rmk:BZFN_monoidal},
gives condition (1) for $\widehat V_0$: $M$ is self-dual over $E$, and
$M\otimes_E M\simeq A\simeq\End_E\left(M\right)$. Finally, $M$ is compactly generated by the
Koszul complex $K_V$ of the origin, so $M\simeq\Ind\left(\Perf_0\left(V\right)\right)$, where
$\Perf_0\left(V\right)$ denotes perfect complexes set-theoretically supported at the origin.
Since $\mathcal{H}om\left(-,\cO_V\right)$ preserves $\Perf_0\left(V\right)$, we have
$\Perf_0\left(V\right)^{\op}\simeq\Perf_0\left(V\right)$ and $M$ is self-dual in $\Pre$.

\emph{Step 2: Condition (2).}
In Lemma \ref{lem:Morita_context}, $\Psi = M\otimes_A-$ is left adjoint to
$\Phi = M^\vee\otimes_E- = \Hom_E\left(M,-\right)$, and it is fully faithful. So the counit
$\varepsilon\colon P=\Psi\Phi\left(E\right)\to E$ is an equivalence as soon as $E$ lies in the
essential image of $\Psi$, which we now show.
Let $\Delta^\times\colon E\to M$ be the right adjoint of $\Delta_*\colon M\to E$. Since $\Delta$ is a
regular closed immersion of codimension $\dim V$, we have
$\Delta^\times\simeq\Delta^*\left(-\right)\otimes\det\left(N_\Delta\right)\left[-\dim V\right]$, so $\Delta^\times$
is continuous and $E$-linear.
It is also conservative, and this is where completing matters: $\Delta^\times$ is not
conservative on $\QC\left(V^2\right)$, but on $E$ it is. Indeed, the compact generator
$K_{V^2}$ of $E$ is $\Delta_*K_V$, since the ideals $\left(x_i,y_i\right)$ and
$\left(x_i-y_i,x_i\right)$ agree. So if $\Delta^\times F = 0$ then
$\Hom\left(K_{V^2},F\right)\simeq\Hom\left(K_V,\Delta^\times F\right)=0$, and $F=0$.
By the Barr--Beck--Lurie theorem \cite[Theorem 4.7.3.5]{L:HA}, $\Delta^\times$ identifies $E$ with
the category $\Mod_T\left(M\right)$ of modules over the monad $T=\Delta^\times\Delta_*$. This monad
is $E$-linear, so it is an algebra in $\End_E\left(M\right)\simeq A$, and
$\Mod_T\left(M\right)\simeq M\otimes_A\Mod_T\left(A\right)$ \cite[Theorem 4.8.4.6]{L:HA}.
Hence $E\simeq\Psi\left(\Mod_T\left(A\right)\right)$ lies in the essential image of $\Psi$.
\end{proof}

\section{Shifted \texorpdfstring{$\cD$}{D}-modules and shifted deformation quantization}
\label{sec:dmod}

The theory of $\cD$-modules, modules over the sheaf of differential operators, is closely
related to the theory of quasi-coherent sheaves in Section \ref{sec:QC}. Of its many treatments
\cite{HTT, GR:v2}, we use the one in which $\cD$-modules on $X$ are ind-coherent sheaves on the
de Rham stack $X_\dR$ \cite{GR:v2}.
This was further generalized to the notion of a \emph{shifted} differential operator (and
correspondingly shifted $\cD$-module) in
\cite{Ber:Sph,Ber:center}, using the sheaf theory $\IndCoh_0$ of Arinkin--Gaitsgory
\cite{AG:crystal}.
Ho and Li \cite{HoLi} generalized this to the $\EE_n$-Hecke categories
$\mathrm{Hecke}_n\left(Y,Z\right)$ attached to a map of stacks $Y\to Z$, which recover
shifted $\cD$-modules for $Z=\pt$, and computed their factorization homology over
arbitrary manifolds.
We will work in this setting.
Using their excision theorem, we compute the higher centers of $\DMod^2\left(X\right)$ as
factorization homology over spheres (Theorem \ref{thm:centers}), alongside the values of
the associated TQFT on closed surfaces \cite{Ber:Sph,HoLi}.
The upshot (Theorem \ref{thm:BW}) is that for $X = V\left[1\right]$, with $V$ a non-zero vector
space, the $\EE_2$-category $\DMod^2\left(X\right)$ is $3$-dualizable and \emph{non-degenerate}
(Definition \ref{defn:inv_E2}) in the Morita $\left(\infty , 4\right)$-category of
$\EE_2$-monoidal categories (see Section \ref{sec:dual_Ek}), but neither factorizable nor
cofactorizable, and hence not invertible there.

\begin{rmk}[Indexing]
\label{rmk:indexing}
Our shift $n$ is the shift of the cotangent bundle being quantized: $\DMod^n\left(X\right)$
is $\EE_n$-monoidal and quantizes $T^*\left[n\right]X$.
In \cite{Ber:Sph} the same category is denoted
$\mathrm{Sph}\left(X , n-1\right)$, and is called the category of
\emph{$n$-shifted $\cD$-modules}.
\end{rmk}

Throughout this section, $k$ is algebraically closed, as in \cite{Ber:ShvCatH}. Following
\cite[Section 2]{Ber:Sph}, a \emph{stack} means a
derived algebraic stack over $k$ which is quasi-compact, locally
of finite presentation, with affine diagonal and perfect cotangent complex.
We also assume stacks are perfect in the sense of \cite{BZFN}, as required in \cite{HoLi};
this is preserved by the mapping stacks below \cite[Remark 3.1.1]{HoLi}.
A map of stacks $Y\to Z$ is \emph{bounded} if $S\times_Z Y$ is bounded for every bounded
affine scheme $S\to Z$ \cite[Section 2.1.5]{Ber:center}, \cite[Section 2.2.3]{HoLi}.

\subsection{Ind-coherent sheaves and formal completions}
\label{sec:IndCoh}

We make a brief aside to define the notion of ind-coherent sheaves, as it will be relevant
throughout this section.

Let $X$ be a derived scheme locally almost of finite type, and write
$\Coh\left(X\right) \subset \QC\left(X\right)$ for the complexes with bounded coherent
cohomology. The category of ind-coherent sheaves on $X$ is defined to be the
ind-completion of this:
\begin{equation*}
\IndCoh\left(X\right)\ceqq \Ind\left(\Coh\left(X\right)\right)
\end{equation*}
This is extended to stacks by taking the limit over affines along $!$-pullback
\cite[Part II]{GR:v1}.
There is a functor
$\Upsilon_X \colon \QC\left(X\right)\to \IndCoh\left(X\right)$,
$F \mapsto F\otimes \omega_X$, intertwining $f^*$ with $f^!$.
It is an equivalence when $X$ is smooth, and fully faithful when $X$ is bounded, i.e.\
eventually coconnective (since then $\omega_X$ is coherent with
$\End\left(\omega_X\right)\simeq\cO_X$). Quasi-smooth stacks are bounded, but taking loops
quickly produces unbounded stacks: $\cL$ sends smooth stacks to quasi-smooth ones, and
quasi-smooth but non-smooth stacks to unbounded ones \cite[Section 2.2]{Ber:Sph}.

Given a map of stacks $f\colon Y\to Z$, not necessarily a closed embedding,
the \emph{formal completion} of $Z$ along $f$ is
$Z^\wedge_Y \ceqq Z\times_{Z_\dR} Y_\dR$,
where $\left(-\right)_\dR$ is the de Rham stack \cite{GR:v2}.
For example, $\pt^\wedge_Y = Y_\dR$, and $Z^\wedge_Y\simeq Z$ whenever $f$ is an
isomorphism on underlying reduced classical stacks; e.g.\
$\left(\cL V\right)^\wedge_V\simeq \cL V$ for a vector space $V$.
Formal completions are not themselves stacks in the above sense. They only enter through
the arrow $Y\to Z$, and $\IndCoh_0$ below is really a sheaf theory on such arrows.

\begin{defn}[{\cite{AG:crystal}, \cite[Section 3]{Ber:center}}]
\label{defn:IndCoh0}
Let $Y\to Z$ be a map of stacks with $Y$ bounded. Then
$\IndCoh_0\left(Z^\wedge_Y\right)$ is defined as the full subcategory of
$\IndCoh\left(Z^\wedge_Y\right)$ given by the fiber product
\begin{equation*}
\IndCoh_0\left(Z^\wedge_Y\right) \ceqq
\IndCoh\left(Z^\wedge_Y\right)\times_{\IndCoh\left(Y\right)} \QC\left(Y\right) \ ,
\end{equation*}
where the maps are $!$-pullback and $\Upsilon_Y$.
If $Y$ is unbounded, this fiber product is not to be considered, and the definition is
modified as in \cite[Definition 4.1.6]{Ber:center}.
\end{defn}

The notation is abusive, as $\IndCoh_0\left(Z^\wedge_Y\right)$ depends on the map $Y\to Z$.
When $Y$ is smooth, $\IndCoh_0\left(Z^\wedge_Y\right) = \IndCoh\left(Z^\wedge_Y\right)$.
$\IndCoh_0$ is functorial for correspondences \cite[Theorem 3.3.3]{Ber:center}, satisfies
K\"unneth \cite[Section 3.4]{Ber:center}, and satisfies descent \cite[Proposition 4.4.1]{Ber:center}, \cite[Proposition 2.3.6]{HoLi};
these are the inputs to the excision theorem \cite[Theorem 3.3.1]{HoLi}.

For $Z = \pt$ we get the \emph{exotic} $\cD$-modules $\IndCoh_0\left(Y_\dR\right)$ of
\cite{Ber:center}. These agree with the usual ones $\IndCoh\left(Y_\dR\right)$ for $Y$
smooth, but unlike the usual ones they can see derived structure.
The only case we will need is the following, in which the derived structure is invisible
after all.

\begin{lem}
\label{lem:DMod0_qsmooth}
Let $Y$ be a bounded stack and $f\colon Y\to Y'$ a map to a smooth stack which
induces an isomorphism $Y_\dR\simeq Y'_\dR$. Then the forgetful functor is an equivalence
\begin{equation*}
\IndCoh_0\left(Y_\dR\right) \lto{\simeq} \IndCoh\left(Y_\dR\right)\simeq
\IndCoh\left(Y'_\dR\right) = \DMod\left(Y'\right) \ .
\end{equation*}
\end{lem}

\begin{proof}
The $!$-pullback $\IndCoh\left(Y_\dR\right)\to\IndCoh\left(Y\right)$ factors as
$\IndCoh\left(Y'_\dR\right)\to\IndCoh\left(Y'\right)\lto{f^!}\IndCoh\left(Y\right)$.
Since $Y'$ is smooth, $\IndCoh\left(Y'\right) = \Upsilon_{Y'}\QC\left(Y'\right)$, and
$f^!\Upsilon_{Y'}\simeq\Upsilon_Y f^*$.
So every object of $\IndCoh\left(Y_\dR\right)$ has underlying ind-coherent sheaf in the
image of $\Upsilon_Y$, which is fully faithful as $Y$ is bounded. Hence the
projection from the fiber product in Definition \ref{defn:IndCoh0} is an equivalence.
\end{proof}

\begin{exm}
\label{exm:DMod0_qsmooth}
Lemma \ref{lem:DMod0_qsmooth} applies to $Y = Y'\times U\left[-1\right]\to Y'$ for
$Y'$ smooth and $U$ a vector space, where $U\left[-1\right]\ceqq\pt\times_U\pt$; and to
$Y = \mathfrak{g}\left[-1\right]/G\to BG$ for an affine algebraic group $G$.
\end{exm}

For a finite-dimensional vector space $U$, write $U\left[1\right]\ceqq BU$ and
$U\left[-1\right]\ceqq\pt\times_U\pt$ for the corresponding \emph{odd vector spaces}. Their
functions are exterior algebras: For $\epsilon=\pm1$,
$\cO\left(U\left[\epsilon\right]\right)\simeq\Lambda^\bullet\left(U^\vee\left[-\epsilon\right]\right)$.

\begin{prop}[Odd vector spaces]
\label{prop:odd}
Let $U$ be a finite-dimensional vector space.
\begin{enumerate}
\item The Weyl algebra of $U\left[\epsilon\right]$, $\epsilon=\pm1$, i.e.\ the algebra generated by the linear
functions $U^\vee\left[-\epsilon\right]$ and the constant vector fields $U\left[\epsilon\right]$, subject to
the graded commutation relations $\left[\xi,x\right]=\langle\xi,x\rangle$, is the Clifford algebra
$\mathrm{Cl}$ of the hyperbolic pairing on $U^\vee\left[-\epsilon\right]\dsum U\left[\epsilon\right]$. Its
action on $\cO\left(U\left[\epsilon\right]\right)$ is an isomorphism
$\mathrm{Cl}\simeq\End_k\left(\Lambda^\bullet\left(U^\vee\left[-\epsilon\right]\right)\right)$, so
$\mathrm{Cl}$ is Morita trivial: $\Mod_{\mathrm{Cl}}\simeq\Vect$.
\item Correspondingly, for $Y$ smooth and $U_1,U_2$ finite-dimensional vector spaces,
\begin{equation*}
\DMod^0\left(Y\times U_1\left[1\right]\times U_2\left[-1\right]\right)\simeq\DMod\left(Y\right) \ .
\end{equation*}
In particular, $\DMod^0$ of an odd vector space is $\Vect$.
\end{enumerate}
\end{prop}

\begin{proof}
For (1): The generators are odd, so graded commutators are anticommutators, and the relations
are $x^2=0$, $\xi^2=0$ and $\xi x+x\xi=\langle\xi,x\rangle$. These are the relations of the Clifford
algebra of the hyperbolic pairing, acting on $\Lambda^\bullet\left(U^\vee\left[-\epsilon\right]\right)$
by multiplication and contraction. This is the Fock representation, which is faithful, and both
sides have dimension $2^{2\dim U}$, so it is an isomorphism. Since
$\Lambda^\bullet\left(U^\vee\left[-\epsilon\right]\right)$ is a non-zero perfect complex, it is a
compact generator of $\Vect$, and $\Mod_{\End\left(P\right)}\simeq\Vect$ for $P$ such a generator.

For (2): $Y\times U_1\left[1\right]$ is smooth, so Lemma \ref{lem:DMod0_qsmooth} (Example
\ref{exm:DMod0_qsmooth}) gives
\begin{equation*}
\DMod^0\left(Y\times U_1\left[1\right]\times U_2\left[-1\right]\right)\simeq
\DMod\left(Y\times U_1\left[1\right]\right) \ .
\end{equation*}
It remains to show that pullback along the projection
$Y\times BU_1\to Y$ is an equivalence. By smooth descent for $\cD$-modules \cite{GR:v2} along
$Y\to Y\times BU_1$, whose \v{C}ech nerve is $Y\times U_1^{\bullet}$, the category
$\DMod\left(Y\times BU_1\right)$ is the limit of the cosimplicial category
$\DMod\left(Y\times U_1^\bullet\right)$. All face maps of this \v{C}ech nerve commute with the
projections to $Y$, and pullback along each projection $Y\times U_1^p\to Y$ is fully faithful, by
the projection formula and $H^*_\dR\left(\AA^m\right)\simeq k$. The component in degree $p$ of an
object of the limit is a pullback of its component in degree $0$, so it lies in the essential
image of $\DMod\left(Y\right)$, and on these essential images the cosimplicial category is
equivalent to the constant one with value $\DMod\left(Y\right)$. Since $\Delta$ is weakly
contractible, the limit is $\DMod\left(Y\right)$.
\end{proof}

\subsection{Shifted \texorpdfstring{$\cD$}{D}-modules}

\begin{defn}[\cite{Ber:Sph}]
\label{defn:DMod_n}
Given a stack $X$ and $n\geq 0$,
the category of $n$-shifted $\cD$-modules on $X$ is
\begin{equation*}
\DMod^n\left(X\right) \ceqq
\IndCoh_0\left(
\left(
X^{S^{n-1}}
\right)^{\wedge}_X
\right) \ ,
\end{equation*}
where $X\to X^{S^{n-1}}$ is the inclusion of constant maps, and $S^{-1}\ceqq\emptyset$.
\end{defn}

The $\EE_{n}$-monoidal structure comes from the fact that this is (sheaves on the
completed) mapping stack from $S^{n-1}$: The multiplication is $\IndCoh_0$-pull-push
along the restriction maps from $X^{P_n}$, where $P_n\ceqq\left(D^n\right)^{\circ\circ}$
is the $n$-dimensional pair of pants, and similarly for more little disks.
For a discussion of the $\mathbb{E}_n$-monoidal structure, we recommend
\cite{Toen:branes}.
$SO\left(n\right)$ acts on all of this data, so $\DMod^n\left(X\right)$ defines an
\emph{oriented} $n$-dimensional TQFT, computed by factorization homology, via the cobordism
hypothesis \cite{L:Cob}, \cite[Section 3.2]{Ber:Sph}. For $n=2$, the only $3$-dimensional theory we
use is the framed one (Lemma \ref{lem:inv_surfaces}).
In the notation of \cite{HoLi}, $\DMod^n\left(X\right)=\mathrm{Hecke}_n\left(X,\pt\right)$.
There the $\EE_n$-structure is constructed as the restriction to disks of a symmetric monoidal
functor out of all $n$-manifolds and open embeddings \cite[Sections 3.1--3.2]{HoLi}, so no
framing is needed.

\begin{exm}
\label{exm:DMod_small_n}
\begin{enumerate}
\item For $n = 0$ we obtain the exotic $\cD$-modules $\DMod^0\left(X\right) =
\IndCoh_0\left(X_\dR\right)$, which are the usual ones when $X$ is smooth.
\item For $n = 1$ we obtain the monoidal category
\begin{equation*}
\DMod^1\left(X\right) = \IndCoh_0\left(\left(X\times X\right)^\wedge_X\right)
\end{equation*}
of \cite{Ber:center}.
\item For $n=2$ we obtain the $\EE_2$-monoidal category
\begin{equation*}
\DMod^2\left(X\right) = \IndCoh_0\left(\left(\cL X\right)^\wedge_X\right).
\end{equation*}
For a vector space $V$ this is $\IndCoh\left(\cL V\right)
= \IndCoh\left(V\times V\left[-1\right]\right)$, and for $X = BG$ it is
$\IndCoh\left(\widehat{G}_e/G\right)$, where $G$ acts on the formal group $\widehat G_e$ by
conjugation.
\end{enumerate}
\end{exm}

\begin{rmk}
\label{rmk:DMod_center}
\begin{enumerate}
\item For $X$ a smooth scheme, $X\to\cL X$ is a nil-isomorphism, so
$\DMod^2\left(X\right)\simeq\IndCoh\left(\cL X\right)$. The convolution category
$\QC\left(\cL X\right)$ of Section \ref{sec:QC} is obtained from it by requiring the singular
support to vanish \cite{AG:singsupp}, i.e.\ by completing along the zero section of
$T^*\left[2\right]X$.
\item Both are Drinfeld centers. $\QC\left(\cL X\right)$ is the center of $\QC\left(X\right)$
(Theorem \ref{thm:BZFN_loops}), i.e.\ of the $2$-category of sheaves of categories on $X$ when
$X$ is $1$-affine. Similarly, $\DMod^2\left(X\right)\simeq\IndCoh\left(\cL X\right)$ is the center of the
$2$-category $2\IndCoh\left(X\right)$ of \cite[Section 1]{Ste:IndCoh}, and it is the value on $S^1$ of
the three-dimensional TQFT built from $2\IndCoh\left(X\right)$ \cite[Corollary 1.3]{Ste:IndCoh}, a
version of Rozansky--Witten theory with target $T^*\left[2\right]X$ (cf.\ \cite[Section D.6]{BZSV}).
\item The symmetric monoidal category $\DMod\left(X\right)$ of ordinary $\cD$-modules acts
centrally on $\DMod^n\left(X\right)$ for every $n\geq1$, so $\DMod^n\left(X\right)$ sheafifies over
$X_\dR$. So invertibility can only be expected when $X$ has no de Rham cohomology, in the sense
that $\DMod\left(X\right)\simeq\Vect$. Proposition \ref{prop:BG}, where
$\cZ_2\left(\DMod^2\left(BG\right)\right)\simeq\DMod\left(BG\right)$, makes this precise for $X=BG$,
and Theorem \ref{thm:BW} shows that invertibility can fail even when $\DMod\left(X\right)\simeq\Vect$.
\end{enumerate}
\end{rmk}

\subsubsection{Factorization homology}

The factorization homology of shifted $\cD$-modules is computed by an excision theorem for
$\IndCoh_0$ \cite[Theorem 3.3.1]{HoLi}, cf.\ \cite[Claim 3.1]{Ber:Sph}.

\begin{prop}[{\cite[Theorem 3.3.1 and Corollary 1.3.12]{HoLi}}]
\label{prop:int_DMod}
Let $M$ be a closed $m$-manifold, $m\leq n$. Then, integrating over $M\times\RR^{n-m}$,
\begin{equation*}
\int_M\DMod^n\left(X\right)\simeq\DMod^{n-m}\left(X^M\right) \ .
\end{equation*}
In particular, $\int_\Sigma\DMod^2\left(X\right)\simeq\DMod^0\left(X^\Sigma\right)=
\IndCoh_0\left(\left(X^\Sigma\right)_\dR\right)$ for a closed oriented surface $\Sigma$.
\end{prop}

This is \cite[Corollary 1.3.12]{HoLi} with $Z=\pt$, since
$\partial\left(M\times D^{n-m}\right)=M\times S^{n-m-1}$ and
$X^{M\times S^{n-m-1}}\simeq\left(X^M\right)^{S^{n-m-1}}$; cf.\ also \cite[Claim 3.1 and
Theorem 5.3]{Ber:Sph}. For $M$ oriented, it is the quantization of the AKSZ
statement \cite{PTVV}
\begin{equation*}
\Maps\left(M , T^*\left[n\right]X\right) \simeq T^*\left[n-m\right]X^M \ .
\end{equation*}

\subsection{Rigidity and \texorpdfstring{$3$}{3}-dualizability}
\label{sec:rigid}

For $n\geq 1$, gluing $S^{n-1}$ from two disks along $S^{n-2}$ gives
$X^{S^{n-1}}\simeq X\times_{X^{S^{n-2}}}X$, so
\begin{equation*}
\DMod^n\left(X\right)\simeq\IndCoh_0\left(\left(X\times_{X^{S^{n-2}}}X\right)^\wedge_X\right) \ ,
\end{equation*}
and the monoidal structure coming from any one of the $n$ directions (these are all equivalent
by $SO\left(n\right)$-invariance) is convolution relative to $X^{S^{n-2}}$ \cite{Ber:center}.
This is the $\EE_1$-Hecke category $\mathrm{Hecke}_1\left(X,X^{S^{n-2}}\right)$ of \cite{HoLi}.
Recall that a monoidal category in $\Pre$ is \emph{rigid} if it is compactly generated, its
unit is compact, and every compact object has a left and a right dual
\cite[Chapter 1, Section 9]{GR:v1}, \cite{G:ShvCat}.
For $n=1$, the following is observed in \cite{Ber:center}; the same argument applies relative
to $X^{S^{n-2}}$.

\begin{prop}
\label{prop:rigid}
Let $X$ be a bounded stack and $n\geq1$. Then
$\DMod^n\left(X\right)$ is rigid, and it is compactly generated by the objects
$\iota\left(P\right)$ for $P\in\Perf\left(X\right)$, where
$\iota\colon\QC\left(X\right)\to\DMod^n\left(X\right)$ is pushforward along the diagonal.
\end{prop}

\begin{proof}
Write $Z\ceqq X^{S^{n-2}}$ and $\Delta\colon X\to\left(X\times_Z X\right)^\wedge_X$ for the diagonal. It induces an
equivalence on de Rham stacks, so it is an inf-schematic nil-isomorphism; in particular it is
ind-proper, so $\Delta_*\dashv\Delta^!$, and $\Delta^!$ is conservative \cite[Part~I]{GR:v2}.
\begin{enumerate}
\item By our conventions $X$ is QCA, so in characteristic zero $\QC\left(X\right)$ is
compactly generated by $\Perf\left(X\right)$, and $\cO_X$ is compact \cite{DG:finiteness}.
\item The functor
\begin{equation*}
\iota = \Delta_*\circ\Upsilon_X
\end{equation*}
lands in $\IndCoh_0$, and is monoidal from $\left(\QC\left(X\right),\otimes\right)$ to the
convolution structure \cite{Ber:center}.
Since $\Upsilon_X$ is fully faithful, $\IndCoh_0$ is the full subcategory of
\begin{equation*}
\IndCoh\left(\left(X\times_ZX\right)^\wedge_X\right)
\end{equation*}
on objects $F$ with
$\Delta^!F\in\Upsilon_X\QC\left(X\right)$. So the right adjoint of $\iota$ is
$\Upsilon_X^{-1}\circ\Delta^!$, the projection in Definition \ref{defn:IndCoh0}, which is
continuous and conservative.
\item Hence $\iota$ preserves compact objects, and the $\iota\left(P\right)$ for
$P\in\Perf\left(X\right)$ form a set of compact generators. In particular the unit
$\iota\left(\cO_X\right)$ is compact.
\item Each $\iota\left(P\right)$ has left and right duals, since $P$ does and $\iota$ is
monoidal. Objects with left (resp.\ right) duals form a thick subcategory, and the compact
objects form the thick subcategory generated by the $\iota\left(P\right)$. So every compact
object has duals. \qedhere
\end{enumerate}
\end{proof}

\begin{cor}
\label{cor:3dual}
Let $X$ be bounded and $\cA = \DMod^2\left(X\right)$. Then $\cA$ and
$\cA^e\simeq\DMod^2\left(X^2\right)$ are rigid.
If moreover $\cL X$ is bounded, e.g.\ if $X$ is smooth, then
$\HC\left(\cA\right)\simeq\DMod^1\left(\cL X\right)$ is rigid,
$\cA$ is dualizable in $\Pre$ and as a module over $\cA^e$ and over $\HC\left(\cA\right)$, and $\cA$
is $3$-dualizable in $\umor_2\left(\Pre\right)$.
\end{cor}

\begin{proof}
Apply Proposition \ref{prop:rigid} to $\DMod^2\left(X\right)$, $\DMod^2\left(X^2\right)$ and
$\DMod^1\left(\cL X\right)$. The identification of
$\HC\left(\cA\right)$ as a monoidal category is \cite[Corollary 5.2]{Ber:Sph}; it is also
\cite[Theorem 3.3.1]{HoLi} for $Z=\pt$ and $M=S^1\times\RR$, with stacking of cylinders
inducing convolution (cf.\ \cite[Example 1.3.8]{HoLi}).
Since $\cA$ is compactly generated, it is dualizable in $\Pre$. A module over a rigid
monoidal category which is dualizable in $\Pre$ is dualizable as a module
\cite[Proposition 2.1.6]{Ber:ShvCatH}, cf.\ \cite[Chapter 1, Section 9]{GR:v1}. Now apply Theorem \ref{thm:inv_E2}.
\end{proof}

\subsection{Higher centers and surfaces}
\label{sec:DMod_centers}

Let $\cA$ be an $\EE_2$-algebra in $\Pre$ which is $SO\left(2\right)$-invariant, so that
factorization homology over oriented surfaces makes sense.

\begin{lem}
\label{lem:inv_surfaces}
If $\cA$ is invertible in $\umor_2\left(\Pre\right)$, then $\int_{T^2}\cA$ is an invertible
object of $\Pre$, for any framing of $T^2$.
\end{lem}

\begin{proof}
The framed TQFT defined by $\cA$ has invertible value on the point, so it factors through the
maximal Picard $\infty$-groupoid of the target, and it takes invertible values on all framed
bordisms.
Its value on the framed torus is $\int_{T^2}\cA$ \cite{Sch:thesis}, an object of
$\Omega^2\umor_2\left(\Pre\right)\simeq\Pre$ (Section \ref{sec:mor}).
\end{proof}

This is a condition on the values $\int_\Sigma\cA$, \emph{not} on $\cA$ itself: An invertible
$\EE_2$-algebra need not be invertible in $\Pre$. For example, a modular tensor category is
invertible in $\umor_2\left(\Pre\right)$ \cite{BJSS} without being equivalent to $\Vect$, but its
factorization homology over every closed surface is $\Vect$ \cite[Theorem 9.7]{KT:4d}
(both in the abelian setting, cf.\ Remark \ref{rmk:abelian}).

\begin{lem}
\label{lem:unit_Vect}
Let $\cZ$ be a monoidal category in $\Pre$ which is equivalent to $\Vect$ as a plain
category. Then the unit map $\Vect\to\cZ$ is an equivalence.
\end{lem}

\begin{proof}
Identify $\cZ$ with $\Vect$. A colimit-preserving $k$-bilinear functor is of the form
\begin{equation*}
\Vect\times\Vect\to\Vect\ ,\quad \left(x,y\right)\mapsto x\otimes y\otimes E
\end{equation*}
with
$E\ceqq k\otimes_\cZ k$. If $\one$ is the unit, then $\one\otimes_\cZ k\simeq k$ says
$\one\otimes E\simeq k$, so $\one$ and $E$ are invertible in $\Vect$, i.e.\ shifts of $k$.
So $x\mapsto x\otimes\one$ is an equivalence.
\end{proof}

The conditions of Definition \ref{defn:inv_E2} concern the higher centers
\begin{equation*}
\cZ_j\left(\cA\right)=\End_{\int_{S^{j-1}\times\RR^{3-j}}\cA}\left(\cA\right)
\end{equation*}
rather than
factorization homology. For shifted $\cD$-modules, $\cZ_j$ is factorization homology over the
$j$-sphere.

\begin{lem}
\label{lem:H_duals}
\leavevmode
\begin{enumerate}
\item For $X$ bounded and $n\geq1$, the monoidal category $\cC=\DMod^n\left(X\right)$ is its own
right dual as a left $\cC^e$-module.
\item Let $f\colon\mathcal{Y}\to\mathcal{Z}$ be a bounded map of bounded stacks. Then the
transfer bimodules between $\DMod^1\left(\mathcal{Y}\right)$ and $\DMod^1\left(\mathcal{Z}\right)$
\begin{equation*}
\mathbb{T}_{\mathcal{Y}\to\mathcal{Z}}\ceqq
\IndCoh_0\left(\left(\mathcal{Y}\times\mathcal{Z}\right)^\wedge_{\mathcal{Y}}\right)
\quad\text{and}\quad
\mathbb{T}_{\mathcal{Z}\leftarrow\mathcal{Y}}\ceqq
\IndCoh_0\left(\left(\mathcal{Z}\times\mathcal{Y}\right)^\wedge_{\mathcal{Y}}\right)
\end{equation*}
are left and right dual to each other \cite{Ber:ShvCatH}.
\end{enumerate}
\end{lem}

\begin{proof}
(1) The category $\cC$ is rigid by Proposition \ref{prop:rigid}. It is also pivotal, since the
proof of \cite[Lemma 3.5.4]{Ber:center} only uses rigidity and the self-duality of
\cite[Proposition 3.2.4]{Ber:center}. So $\cC$ is very rigid \cite[Section 5.1.4]{Ber:center},
i.e.\ the self-duality $\cC\simeq\cC^\vee$ is $\left(\cC,\cC\right)$-bilinear.
By \cite[Corollary 2.1.7]{Ber:ShvCatH}, the right dual of $\cC$ as a module over the rigid
category $\cC^e$ is $\cC^\vee$, so it is $\cC$.

(2) These are the values of the $2$-functor of \cite[Theorem 6.5.3]{Ber:ShvCatH} on the correspondences
$\left[\mathcal{Y}=\mathcal{Y}\to\mathcal{Z}\right]$ and $\left[\mathcal{Z}\leftarrow\mathcal{Y}=\mathcal{Y}\right]$.
Tensoring with them realizes pullback $f^*$ and pushforward $f_*$ between
$\DMod^1\left(\mathcal{Z}\right)$- and $\DMod^1\left(\mathcal{Y}\right)$-modules
\cite[Corollary 6.5.4]{Ber:ShvCatH}, and $f_!\dashv f^*\dashv f_*$ with
$f_!\simeq f_*$ \cite[Corollary 6.5.5]{Ber:ShvCatH}.
\end{proof}

\begin{theorem}
\label{thm:centers}
Let $X$ be smooth and $\cA = \DMod^2\left(X\right)$. For $j=0,1$, and for $j=2$ if the constant
loops $X\to\cL X$ are bounded,\footnote{This holds e.g.\ for $X=BG$ with $G$ an affine algebraic
group, since the pullback of $BG\to G/G$ along the smooth cover $G\to G/G$ is $e\to G$. It fails
for smooth schemes of positive dimension, for which $X^{S^2}\simeq T\left[-2\right]X$ is unbounded
(cf.\ \cite[Section 1.3.13]{HoLi}). We expect this assumption to be unnecessary.}
there are equivalences in $\Pre$
\begin{equation*}
\cZ_j\left(\cA\right)\simeq\int_{S^j\times\RR^{2-j}}\cA\simeq\DMod^{2-j}\left(X^{S^j}\right) \ .
\end{equation*}
For $j=0$ the middle term is $\cA\otimes\cA$, and for $j=1$ it is $\HC\left(\cA\right)$.
\end{theorem}

\begin{proof}
The second equivalence is Proposition \ref{prop:int_DMod}. For the first equivalence, let
$\cB_j\ceqq\int_{S^{j-1}\times\RR^{3-j}}\cA$, so that $\cZ_j\left(\cA\right)=\End_{\cB_j}\left(\cA\right)$.
If $\cA$ has a right dual $\cA^R$ as a left $\cB_j$-module, then
\begin{equation*}
\End_{\cB_j}\left(\cA\right)\simeq\cA^R\otimes_{\cB_j}\cA \ ,
\end{equation*}
so it suffices to find $\cA^R$.

For $j=0$, $\cA^R=\cA^\vee\simeq\cA$ \cite[Proposition 3.2.4]{Ber:center}, so
$\cZ_0\left(\cA\right)\simeq\cA\otimes\cA$.

For $j=1$, $\cB_1=\cA^e$, so $\cA^R\simeq\cA$ by Lemma
\ref{lem:H_duals}(1), and $\cZ_1\left(\cA\right)\simeq\cA\otimes_{\cA^e}\cA=\HC\left(\cA\right)$.

For $j=2$, $\cB_2=\HC\left(\cA\right)\simeq\DMod^1\left(\cL X\right)$ (Corollary \ref{cor:3dual}).
Gluing the collar $S^1\times\RR$ onto the disk gives the correspondence
$\cL X\leftarrow X\to\pt$ of mapping stacks \cite[Section 3.1]{HoLi}, so $\cA$ is the value of
the $2$-functor of \cite[Theorem 6.5.3]{Ber:ShvCatH} on it. Factoring through $X$,
\begin{equation*}
\cA\simeq\mathbb{T}_{\cL X\leftarrow X}\otimes_{\DMod^1\left(X\right)}\mathbb{T}_{X\to\pt} \ .
\end{equation*}
By Lemma \ref{lem:H_duals}(2), applied to $X\to\cL X$ and $X\to\pt$, the right dual $\cA^R$ is the
value on the reversed correspondence $\pt\leftarrow X\to\cL X$. So $\cZ_2\left(\cA\right)$ is the
value on the composite $\pt\leftarrow X\times_{\cL X}X\to\pt$, which is
$\DMod^0\left(X^{S^2}\right)$.
\end{proof}

\begin{rmk}
\label{rmk:centers_spheres}
By Theorem \ref{thm:centers}, for $X$ as there the comparison maps of Corollary
\ref{cor:SSS_inv} for $\cA=\DMod^2\left(X\right)$ are, up to the equivalences of that theorem
(which are not claimed to be compatible with these maps), maps between factorization homologies
over spheres,
\begin{equation*}
\int_{S^{i-1}\times\RR^{3-i}}\cA\longrightarrow\int_{S^{2-i}\times\RR^{i}}\cA \ .
\end{equation*}
For $i=0$ this is $\Vect\to\int_{S^2}\cA$ (non-degeneracy), for $i=1$ it is
$\cA\otimes\cA^{\sigma\rev}\simeq\int_{S^0\times\RR^2}\cA\to\int_{S^1\times\RR}\cA$
(factorizability), and for $i=2$ it is $\int_{S^1\times\RR}\cA\to\int_{S^0\times\RR^2}\cA$
(cofactorizability). So factorizability and cofactorizability compare the same two categories,
$\DMod^2\left(X^2\right)$ and $\DMod^1\left(\cL X\right)$, in opposite directions.
\end{rmk}

\subsection{The critical case: \texorpdfstring{$X = V\left[1\right]$}{X = V[1]}}

Let $V$ be a finite-dimensional vector space, and set $X\ceqq V\left[1\right]$ and
$\cA \ceqq \DMod^2\left(X\right)$.

\begin{theorem}
\label{thm:BW}
For a closed oriented surface $\Sigma_g$ of genus $g$,
\begin{equation*}
\int_{\Sigma_g}\DMod^2\left(V\left[1\right]\right)\simeq\DMod\left(V^{2g}\right) \ .
\end{equation*}
Consequently:
\begin{enumerate}
\item $\cA$ is $3$-dualizable and non-degenerate: $\cZ_2\left(\cA\right)\simeq\Vect$;
\item if $V\neq 0$ then $\int_{T^2}\cA$ is not invertible, so $\cA$ is not invertible;
\item if $V\neq 0$ then $\cA$ is neither factorizable nor cofactorizable.
\end{enumerate}
\end{theorem}

\begin{proof}
Since $V\left[1\right]$ is a linear stack, we have
$X^{\Sigma_g}\simeq V\otimes C^*\left(\Sigma_g\right)\left[1\right]$, i.e.
\begin{equation*}
\left(V\left[1\right]\right)^{\Sigma_g}\simeq
V\left[1\right]\dsum V^{2g}\dsum V\left[-1\right]
\simeq V\left[1\right]\times V^{2g}\times V\left[-1\right] \ .
\end{equation*}
This is quasi-smooth, hence bounded, and the projection to the smooth stack $V\left[1\right]\times V^{2g}$ is an
isomorphism on de Rham stacks. The other stacks involved, $V\left[1\right]$,
$\cL\left(V\left[1\right]\right)\simeq V\times V\left[1\right]$ and
$V\left[1\right]\times V\left[1\right]$, are smooth. So every stack that appears is bounded, and Definition
\ref{defn:IndCoh0} applies as stated. (In general, this is the Hecke pair condition of
\cite[Definition 3.4.2]{HoLi}, which holds for $BH\to BG$ for any homomorphism of affine
algebraic groups \cite[Lemma 3.4.4]{HoLi}, here $BV\to\pt$.) By Proposition \ref{prop:int_DMod} and Proposition \ref{prop:odd}(2), with $Y=V^{2g}$ and
$U_1=U_2=V$,
\begin{equation*}
\int_{\Sigma_g}\cA\simeq\DMod^0\left(\left(V\left[1\right]\right)^{\Sigma_g}\right)
\simeq\DMod\left(V^{2g}\right) \ .
\end{equation*}

For (1): $V\left[1\right]$ is smooth, so $\cA$ is $3$-dualizable by Corollary \ref{cor:3dual}.
Since $V\left[1\right]=BV$, the constant loops $V\left[1\right]\to\cL\left(V\left[1\right]\right)$
are bounded (see the footnote to Theorem \ref{thm:centers}), so
$\cZ_2\left(\cA\right)\simeq\int_{S^2}\cA\simeq\Vect$ by Theorem \ref{thm:centers} and the
case $g=0$. The unit map is then an equivalence by Lemma \ref{lem:unit_Vect}.
For (2) take $g=1$: The value of the TQFT on the torus is
$\int_{T^2}\cA\simeq\DMod\left(V^2\right)\simeq\Mod_D$, with $D$ the Weyl algebra
of $V^2$. This object of $\Pre$ is not invertible: Otherwise
$\Mod_{D^\op\otimes D}\simeq\Mod_D^\vee\otimes\Mod_D\simeq\Vect$, but the compact object
$D^\op\otimes D$ has infinite-dimensional endomorphisms. By Lemma
\ref{lem:inv_surfaces}, $\cA$ is not invertible.

For (3), call a compactly generated category \emph{proper} if all mapping complexes between
compact objects have finite-dimensional total cohomology. This is invariant under equivalence in
$\Pre$, and it holds if there is a compact generator with finite-dimensional endomorphisms. By
Remark \ref{rmk:centers_spheres}, it suffices to show that $\int_{S^0\times\RR^2}\cA$ is proper
and $\int_{S^1\times\RR}\cA$ is not.
\begin{itemize}
\item $\int_{S^0\times\RR^2}\cA\simeq\DMod^2\left(X^2\right)$ is proper. Write $U\ceqq V^2$, so
$X^2\simeq BU$ and $\left(\cL BU\right)^\wedge_{BU}\simeq\widehat U_0\times BU$, with $\widehat U_0$
the formal completion of $U$ at $0$ (Example \ref{exm:DMod_small_n}), so
$\DMod^2\left(X^2\right)\simeq\IndCoh\left(\widehat U_0\right)\otimes\QC\left(BU\right)$. Since $U$
is unipotent, $\Perf\left(BU\right)$ is generated by $\cO_{BU}$, so by Proposition
\ref{prop:rigid} this is compactly generated by
$\iota\left(\cO_{BU}\right)=k_0\boxtimes\cO_{BU}$. Its endomorphisms are
$\End\left(k_0\right)\otimes C^*\left(BU\right)\simeq
\Lambda^\bullet\left(U\left[-1\right]\right)\otimes\Lambda^\bullet\left(U^\vee\left[-1\right]\right)$,
which is finite-dimensional.
\item $\int_{S^1\times\RR}\cA\simeq\DMod^1\left(\cL X\right)$ is
not proper. Its unit $\iota\left(\cO_{\cL X}\right)$ is compact by Proposition \ref{prop:rigid}.
The action of $\DMod^1\left(\cL X\right)$ on $\IndCoh\left(\cL X\right)$ by integral
transforms restricts along $\iota$ to the action of $\QC\left(\cL X\right)$ by tensoring, so the
composite
$\Gamma\left(\cL X,\cO\right)\to\End\left(\iota\left(\cO_{\cL X}\right)\right)\to
\End\left(\omega_{\cL X}\right)\simeq\Gamma\left(\cL X,\cO\right)$ is the identity. Hence
$\Gamma\left(\cL X,\cO\right)\simeq\cO\left(V\right)\otimes\Lambda^\bullet\left(V^\vee\left[-1\right]\right)$,
which is infinite-dimensional as $V\neq0$, is a retract of the endomorphisms of a compact object.
\end{itemize}
So neither the factorizability map
$\cA\otimes\cA^{\sigma\rev}\simeq\int_{S^0\times\RR^2}\cA\to\cZ_1\left(\cA\right)\simeq\int_{S^1\times\RR}\cA$
nor the cofactorizability map
$\HC\left(\cA\right)=\int_{S^1\times\RR}\cA\to\cZ_0\left(\cA\right)\simeq\int_{S^0\times\RR^2}\cA$
is an equivalence.
\end{proof}

Theorem \ref{thm:BW} is an instance of a parity phenomenon. For a linear
stack $X=U$, i.e.\ a graded vector space, $X^M\simeq U\otimes H^*\left(M\right)$. By
Proposition \ref{prop:odd}, odd directions contribute Clifford algebras, which are Morita trivial,
while even directions contribute Weyl algebras, which are not. Since $H^*\left(S^0\right)$ and
$H^*\left(S^2\right)$ are concentrated in even degrees, $U\otimes H^*\left(S^0\right)$ and
$U\otimes H^*\left(S^2\right)$ are odd when $U$ is, so the $S^2$ test of non-degeneracy is passed.
But $H^*\left(S^1\right)$ and $H^*\left(T^2\right)$ have odd classes, so $U\otimes H^*\left(S^1\right)$
and $U\otimes H^*\left(T^2\right)$ have non-zero even part unless $U=0$. For $U=V\left[1\right]$,
this even part is the factor $V$ of $\cL X\simeq V\times V\left[1\right]$, which makes
$\int_{S^1\times\RR}\cA$ non-proper and so breaks both factorizability and cofactorizability, and
the factor $V^{2g}$ of
$\left(V\left[1\right]\right)^{\Sigma_g}\simeq V\left[1\right]\times V^{2g}\times V\left[-1\right]$,
which makes the $T^2$ test fail.

\begin{prop}
\label{prop:BG}
Let $G$ be an affine algebraic group. Then
$\cZ_2\left(\DMod^2\left(BG\right)\right)\simeq\DMod\left(BG\right)$. Hence
$\DMod^2\left(BG\right)$ is non-degenerate if and only if $G$ is unipotent.
\end{prop}

\begin{proof}
Since $BG$ is smooth and $BG\to\cL BG\simeq G/G$ is bounded, Theorem \ref{thm:centers} applies. With
$BG^{S^2}\simeq\mathfrak{g}\left[-1\right]/G$ and Example \ref{exm:DMod0_qsmooth}, it gives
\begin{equation*}
\cZ_2\left(\DMod^2\left(BG\right)\right)\simeq\DMod^0\left(\mathfrak{g}\left[-1\right]/G\right)
\simeq\DMod\left(BG\right) \ .
\end{equation*}
By Lemma \ref{lem:unit_Vect}, it follows that $\DMod^2\left(BG\right)$ is non-degenerate
exactly when $\DMod\left(BG\right)\simeq\Vect$.
If $G$ is unipotent, then $\DMod\left(BG\right)\simeq\Vect$ as in the proof of Proposition
\ref{prop:odd}(2), since $G\simeq\AA^{\dim G}$ as a variety.
Conversely, suppose $\DMod\left(BG\right)\simeq\Vect$.
We may assume $k=\CC$: $G$ is defined over a finitely generated subfield, which embeds into $\CC$,
unipotence is invariant under field extension, and so is the condition that the unit
$\Vect\to\DMod\left(BG\right)$ is an equivalence, since $\DMod$ commutes with extension of scalars and
extension of scalars is conservative.
Let $K$ be a maximal compact subgroup of $G\left(\CC\right)$. The endomorphisms of the unit are
the de Rham cohomology of $BG$, which by the comparison theorem is
$H^*\left(BG\left(\CC\right);\CC\right)\simeq H^*\left(BK;\CC\right)$, since $G\left(\CC\right)$
deformation retracts onto $K$. If $K^\circ\neq1$, a $K$-invariant inner product on its Lie
algebra gives a non-zero class in $H^4\left(BK;\CC\right)$, so $K$ is finite. Then the identity
component $G^\circ$ is unipotent and
$\Gamma\ceqq\pi_0\left(G\right)\simeq\pi_0\left(G\left(\CC\right)\right)\simeq K$. Since
$\DMod\left(BG^\circ\right)\simeq\Vect$, pullback along $q\colon BG\to B\Gamma$ is fully faithful: By
the projection formula and base change along $\pt\to B\Gamma$, the unit $\Id\to q_*q^*$ is
tensoring with the de Rham cohomology of the fiber $BG^\circ$, which is $k$. So
$\Rep\left(\Gamma\right)=\DMod\left(B\Gamma\right)$ embeds fully faithfully into $\Vect$. In
$\Rep\left(\Gamma\right)$ distinct irreducibles have zero mapping complexes, whereas any two
non-zero objects of $\Vect$ have a non-zero mapping complex, so $\Gamma$ is trivial and $G$ is
unipotent.
\end{proof}

For example, $\DMod\left(B\Gamma\right)=\Rep\left(\Gamma\right)$ for $\Gamma$ finite,
$\DMod\left(BT\right)$ is modules over an exterior algebra for a torus $T$, and in general
$\DMod\left(BG\right)\simeq\DMod\left(BG^{\mathrm{red}}\right)$ for $G^{\mathrm{red}}\ceqq G/R_u\left(G\right)$
the reductive quotient.

\begin{cor}
\label{cor:BG_not_inv}
Let $G$ be a non-trivial affine algebraic group which is either not unipotent or abelian.
Then $\DMod^2\left(BG\right)$ is not invertible.
\end{cor}

\begin{proof}
If $G$ is not unipotent, then $\DMod^2\left(BG\right)$ is degenerate by Proposition
\ref{prop:BG}, so not invertible by Theorem \ref{thm:inv_E2}. If $G$ is abelian unipotent, then
$G$ is the additive group of a vector space $V$, so $BG\simeq V\left[1\right]$, and this is Theorem \ref{thm:BW}(2).
\end{proof}

\section{Gauging and anomalies}
\label{sec:gauging}

Throughout this section, $G$ is an affine algebraic group of finite type over a field of
characteristic zero, and $\Rep\left(G\right) = \QC\left(BG\right)$.

Sections \ref{sec:shvcat}--\ref{sec:gauging_anomaly} work in complete generality: $A$ is an
arbitrary $\EE_2$-algebra in $\Pre$, and $G$ acts on it by $\EE_2$-functors. We explain how such
an action produces an object $A^G$ of $\umor_2\left(\Mod_{\Rep\left(G\right)}\left(\Pre\right)\right)$,
and what this gives when $A$ is invertible. In Section \ref{sec:gauging_lines} we specialize to the
case where $A$ is the category of line operators of a $3$-dimensional TQFT, and show that actions
induced by an action on the underlying stack are anomaly-free. In Section \ref{sec:tori} we show that the categories of representations
of quantum groups of tori at roots of unity are invertible relative to quantum Frobenius.
Finally, in Section \ref{sec:aggregated} we discuss an analytic setting.
We begin by recalling sheaves of categories and $1$-affineness in Section \ref{sec:shvcat}.

\subsection{Sheaves of categories and \texorpdfstring{$1$}{1}-affineness}
\label{sec:shvcat}

An ordinary quasi-coherent sheaf on $X$ can be thought of as a functorial assignment
(for every affine $\Spec \left(A\right)$ over $X$) of an $A$-module.
Similarly, a sheaf of categories is a functorial assignment, for every affine $S \to X$,
of a certain module category (in $\Pre$) over the symmetric-monoidal category
$\QC\left(S\right)$.

Formally, $\ShvCat\left(X\right)$ is the limit of $\Mod_{\QC\left(S\right)}\left(\Pre\right)$ over affines
$S\to X$, and $X$ is \emph{$1$-affine} if taking global sections
$\ShvCat\left(X\right)\to\Mod_{\QC\left(X\right)}\left(\Pre\right)$ is an equivalence \cite{G:ShvCat}.
For example, $BG$ is $1$-affine \cite[Theorem 2.2.2]{G:ShvCat}:
$\ShvCat\left(BG\right)\simeq\Mod_{\Rep\left(G\right)}\left(\Pre\right)$.

\subsection{Equivariantization}
\label{sec:gauging_general}

Let $A\in\umor_2\left(\Pre\right)$ be arbitrary.
A sheaf of categories on $BG$ is a $\QC\left(S\right)$-linear category for every
$S\to BG$, compatible with base change.
Computing on the \v{C}ech nerve of $\pt\to BG$, this amounts to a category with a coaction of
$\QC\left(G\right)$.

\begin{defn}
An \emph{action of $G$ on $A$ by $\EE_2$-functors} is an $\EE_2$-algebra $A/G$ in
$\ShvCat\left(BG\right)$ with fiber $A$ at $\pt\to BG$.
Equivalently, it is an $\EE_2$-monoidal coaction
$\rho\colon A\to\QC\left(G\right)\otimes A$, coassociative up to coherent homotopy.
\end{defn}

We write $A/G$ for the sheaf of categories, by analogy with the quotient stack: For a
$G$-stack $X$, the action of $G$ on $\QC\left(X\right)$ is the sheaf of categories on $BG$
whose global sections are $\QC\left(X/G\right)$.
An action is a family of autoequivalences parametrized by the scheme $G$, not merely a
homomorphism from the discrete group $G\left(k\right)$.

\begin{defn}
The \emph{equivariantization} of the action is $A^G\ceqq\Gamma\left(BG,A/G\right)$.
This is the category of objects $a\in A$ equipped with an equivalence
$\rho\left(a\right)\simeq\cO_G\boxtimes a$ satisfying the cocycle condition.
\end{defn}

The functor $\Gamma\left(BG,-\right)\colon\ShvCat\left(BG\right)\to\Mod_{\Rep\left(G\right)}\left(\Pre\right)$ is right
adjoint to the symmetric monoidal functor
$M\mapsto\left(\QC\left(S\right)\otimes_{\Rep\left(G\right)}M\right)_{S\to BG}$.
It is therefore lax symmetric monoidal, and so it sends $\EE_2$-algebras to $\EE_2$-algebras.
Hence $A^G\in\umor_2\left(\Mod_{\Rep\left(G\right)}\left(\Pre\right)\right)$.
Concretely, equivariant objects are multiplied and braided as in $A$.
This is possible because $\rho$ is braided monoidal.
The unit $\Rep\left(G\right)\to A^G$ sends $V$ to $\one_A$ with equivariant structure $V$.
These objects braid trivially with everything, recovering
$\Rep\left(G\right)\subseteq \mathcal{Z}_2\left(A^G\right)$ from the finite setting \cite{DGNO}.

\begin{prop}
\label{prop:equiv}
\leavevmode
\begin{enumerate}
\item Equivariantization defines an equivalence
\begin{equation*}
\umor_2\left(\ShvCat\left(BG\right)\right) \lto{\sim}
\umor_2\left(\Mod_{\Rep\left(G\right)}\left(\Pre\right)\right)\ .
\end{equation*}
Its inverse sends $A'$ to
$\QC\left(\pt\right)\otimes_{\Rep\left(G\right)}A'$ with its induced coaction.
\item $A$ is invertible in $\umor_2\left(\Pre\right)$ if and only if $A^G$ is invertible in
\begin{equation*}
\umor_2\left(\Mod_{\Rep\left(G\right)}\left(\Pre\right)\right).
\end{equation*}
\end{enumerate}
\end{prop}

\begin{proof}
For (1): By $1$-affineness of $BG$ \cite[Theorem 2.2.2]{G:ShvCat},
$\Gamma\left(BG,-\right)$ is a symmetric monoidal equivalence.
The inverse is identified using
\begin{equation*}
\QC\left(\pt\right)\otimes_{\Rep\left(G\right)}\QC\left(\pt\right)\simeq\QC\left(G\right)
\end{equation*}
(Theorem \ref{thm:BZFN}).
This part holds for every $\EE_2$-algebra $A$.

For (2), first assume $A$ is invertible.
By (1), it suffices to show that $A/G$ is invertible in
$\umor_2\left(\ShvCat\left(BG\right)\right)$.
We check the criterion of Theorem \ref{thm:inv_E2} in
$\coef = \ShvCat\left(BG\right)$, using descent: By \cite[Theorem 1.5.7]{G:ShvCat},
\begin{equation*}
\ShvCat\left(BG\right)\simeq\mathrm{Tot}\left(\Mod_{\QC\left(G^\bullet\right)}\left(\Pre\right)\right)
\end{equation*}
as symmetric monoidal $\infty$-categories, with transition functors given by base change,
which are symmetric monoidal and preserve colimits. Write $p_n^*$ for the projections. They
commute with tensor products, relative tensor products and factorization homology, so the
images under $p_n^*$ of $A/G$, $\left(A/G\right)^e$, $\HC\left(A/G\right)$ and their modules
are the base changes $\QC\left(G^n\right)\otimes\left(-\right)$ of the corresponding objects for
$A$.
\begin{itemize}
\item \emph{$3$-dualizability.} An object of a limit of symmetric monoidal $\infty$-categories
is dualizable if and only if all its images are \cite[Lemma 1.4.6]{G:ShvCat}, since the space of
duality data is the limit of the spaces of duality data. The same argument applies to
dualizability of modules. Base change preserves duals, so the dualizability conditions for
$A$ imply those for $A/G$.
\item \emph{The three maps.} Once the relevant modules are dualizable, the centers
$\mathcal{Z}_j\left(A/G\right)$ are computed by duals and relative tensor products, hence commute with
$p_n^*$. A map in a limit of $\infty$-categories is an equivalence if and only if all its images
are, and the images are base changes of the corresponding equivalences for $A$.
\end{itemize}
Conversely, if $A^G$ is invertible, then so is $A/G$ by (1). Taking the fiber at
$\pt\to BG$ is a functor from $\ShvCat\left(BG\right)$ to $\Pre$ which is symmetric monoidal and preserves colimits, so it
induces a symmetric monoidal functor on Morita categories, and sends the invertible $A/G$ to $A$.
\end{proof}

\begin{rmk}
\label{rmk:inv_descent}
The proof of (2) shows in general that an $\EE_2$-algebra in $\ShvCat\left(BG\right)$ is
invertible if and only if its fiber at $\pt\to BG$ is invertible: All conditions of Theorem
\ref{thm:inv_E2} can be checked after pullback along $\pt\to BG$.
\end{rmk}

\subsection{The anomaly of an action on an invertible algebra}
\label{sec:gauging_anomaly}

We remain in the general setting, but now assume that $A$ is invertible, so that $A^G$ is
invertible by Proposition \ref{prop:equiv}(2). We write $\alpha_A$ and $\alpha_{A^G}$ for the
corresponding invertible $4$-dimensional theories \cite{BJSS}.

\begin{nota}
\label{nota:Pic}
For a symmetric monoidal higher category $\coef$, we write $\mathrm{Pic}\left(\coef\right)$ for its
Picard $\infty$-groupoid: the invertible objects of $\coef$, and the invertible morphisms between
them.
\end{nota}

\begin{defn}
\label{defn:Witt}
Let $\coef$ be a symmetric monoidal $\infty$-category. The \emph{Witt group} of $\coef$ is the
group
\begin{equation*}
\mathrm{W}\left(\coef\right)\ceqq\pi_0\mathrm{Pic}\left(\umor_2\left(\coef\right)\right)
\end{equation*}
of Morita classes of invertible $\EE_2$-algebras in $\coef$, under $\otimes$. We use this for
$\coef$ equal to $\Pre$, $\ShvCat\left(BG\right)$ and $\Mod_{\Rep\left(G\right)}\left(\Pre\right)$.
\end{defn}

The inverse of $A$ is $A^{\sigma\rev}$, so we may form the invertible
$\Rep\left(G\right)$-linear $\EE_2$-algebra
\begin{equation*}
\omega\left(A/G\right)\ceqq
A^G\otimes_{\Rep\left(G\right)}\left(\Rep\left(G\right)\otimes A^{\sigma\rev}\right) \ ,
\end{equation*}
which compares $A^G$ with $A$ carrying the trivial action.
Its fiber $A\otimes A^{\sigma\rev}\simeq\one$ is Morita trivial.

\begin{defn}
\label{defn:anomaly}
The \emph{anomaly} of the action is the class
\begin{equation*}
\left[\omega\left(A/G\right)\right]
=\left[A^G\right]-\left[\Rep\left(G\right)\otimes A\right]
\in\mathrm{W}\left(\Mod_{\Rep\left(G\right)}\left(\Pre\right)\right)
\simeq\mathrm{W}\left(\ShvCat\left(BG\right)\right) \ ,
\end{equation*}
where the equivalence is Proposition \ref{prop:equiv}(1).
\end{defn}

The anomaly vanishes exactly when $A^G$ is Morita equivalent over $\Rep\left(G\right)$ to
$\Rep\left(G\right)\otimes A$.

\subsection{Line operators of a \texorpdfstring{$3$}{3}-dimensional theory}
\label{sec:gauging_lines}

We now specialize to line operators (Section \ref{sec:line_ops}).
Let $\mathcal{F}$ be a $3$-dimensional TQFT valued in $\umor_1\left(\Pre\right)$, with $B\ceqq \mathcal{F}\left(\pt\right)$
and category of line operators $A\ceqq \mathcal{F}\left(S^1\right)\simeq \mathcal{Z}_1\left(B\right)$.

A \emph{$G$-symmetry} of $\mathcal{F}$ is a lift of $\mathcal{F}$ to a TQFT $\mathcal{F}_{BG}$ valued in
$\umor_1\left(\ShvCat\left(BG\right)\right)$, i.e.\ a family of theories over $BG$.
On the circle, $\mathcal{F}_{BG}\left(S^1\right)$ is an action of $G$ on $A$ by $\EE_2$-functors,
so the constructions of Sections \ref{sec:gauging_general} and \ref{sec:gauging_anomaly}
apply.
Composing with the equivalence $\Gamma\left(BG,-\right)$ gives a $\Rep\left(G\right)$-linear
TQFT $\mathcal{F}^G$ with
\begin{equation*}
\mathcal{F}^G\left(\pt\right) = B^G \ ,
\qquad
\mathcal{F}^G\left(S^1\right) = A^G \ ,
\end{equation*}
and $B^G$ is a $1$-morphism $A^G\to\one$ in $\umor_2\left(\Mod_{\Rep\left(G\right)}\left(\Pre\right)\right)$.
When $A$ is invertible, $\mathcal{F}^G$ is a boundary condition for $\alpha_{A^G}$.
For example, an action of $G$ on $X$ gives $\QC\left(X\right)^G = \QC\left(X/G\right)$
and $\QC\left(\cL X\right)^G = \QC\left(\cL X/G\right)$.
Here $\cL X/G$ is the relative loop stack of $X/G\to BG$, with relative convolution.

\begin{lem}[Geometric actions are anomaly-free]
\label{lem:geometric}
Let $G$ act on a perfect stack $X$ satisfying diagonal descent (Definition
\ref{defn:diag_descent}), with $X/G$ perfect, so
that $A=\left(\QC\left(\cL X\right),*\right)$ is Morita trivial (Theorem \ref{thm:trivial_QC}). Then
$\QC\left(X/G\right)$ is an invertible $1$-morphism from
$\QC\left(\cL X/G\right)\simeq A^G$ to $\Rep\left(G\right)$ in $\umor_2\left(\Mod_{\Rep\left(G\right)}\left(\Pre\right)\right)$,
and the anomaly of the action of $G$ on $A$ vanishes.
\end{lem}

\begin{proof}
Write $\cL X/G\simeq\left(X/G\right)\times_{\left(X\times X\right)/G}\left(X/G\right)$ for the
relative loop stack of $X/G\to BG$. Theorem \ref{thm:BZFN} for the diagonal of $X/G$ over $BG$
identifies the Drinfeld center of the commutative algebra $\QC\left(X/G\right)$ in
$\Mod_{\Rep\left(G\right)}\left(\Pre\right)$ with $\QC\left(\cL X/G\right)\simeq A^G$. So
$\QC\left(X/G\right)$ is a $1$-morphism $A^G\to\Rep\left(G\right)$ in
$\umor_2\left(\Mod_{\Rep\left(G\right)}\left(\Pre\right)\right)$. Its fiber at $\pt\to BG$ is
$\QC\left(X\right)\colon\QC\left(\cL X\right)\to\one$, since
$\QC\left(\pt\right)\otimes_{\Rep\left(G\right)}\QC\left(Y/G\right)\simeq\QC\left(Y\right)$ for
$Y=X,\cL X$ (Theorem \ref{thm:BZFN}). This fiber is invertible by Theorem \ref{thm:trivial_QC},
and invertibility can be checked after pullback along $\pt\to BG$: As in the proof of
Proposition \ref{prop:equiv}(2), the dualizability conditions and the four maps of Theorem
\ref{thm:trivial_Alg2} are computed termwise in
$\ShvCat\left(BG\right)\simeq\mathrm{Tot}\left(\Mod_{\QC\left(G^\bullet\right)}\left(\Pre\right)\right)$,
and a map in a limit is an equivalence if all its images are. So $\left[A^G\right]=\left[\Rep\left(G\right)\right]=\left[\Rep\left(G\right)\otimes A\right]$,
and $\left[\omega\left(A/G\right)\right]=0$.
\end{proof}

In terms of the $2$-shifted cotangent bundle $T^*\left[2\right]X$ (cf.\ Remark
\ref{rmk:DMod_center}), Lemma \ref{lem:geometric} says that symmetries of $T^*\left[2\right]X$
induced by symmetries of $X$ are anomaly-free.

\subsection{Quantum groups of tori at roots of unity}
\label{sec:detect}

We now apply the constructions of this section to quantum groups of tori at roots of unity, where
the quantum Frobenius gives a central $\EE_2$-functor from $\Rep\left(H\right)$, for a torus $H$, to
the category of representations, landing in transparent objects. We show that the category of
representations is invertible over $\Rep\left(H\right)$ (Theorem \ref{thm:T_invertible}), so that
the anomaly of the action is defined (Definition \ref{defn:anomaly}).

Throughout this subsection, $k=\CC$, and the group acting by $\EE_2$-functors is denoted $H$.

\label{sec:tori}

Let $T$ be a torus with character lattice $\Lambda$, and let $Q\colon\Lambda\to k^\times$ be a
quadratic form, with double braiding
$b\left(\lambda,\mu\right)=Q\left(\lambda+\mu\right)Q\left(\lambda\right)^{-1}Q\left(\mu\right)^{-1}$
and radical $\Lambda_0$ of $b$. Assume that $Q|_{\Lambda_0}=1$ and that $\Lambda/\Lambda_0$ is
finite, and let $H\ceqq\Hom\left(\Lambda_0,\GG_m\right)$. By \cite{JS:braided}, $Q$ determines a
braided structure on $\Lambda$-graded vector spaces, and we write $\Rep\left(T_q\right)$ for its
derived category, i.e.\ $\Lambda$-graded complexes. The basic example is $T=\GG_m$,
$Q\left(a\right)=q^{a^2}$ with $q$ a primitive $n$-th root of unity, $n$ odd, for which
$\Lambda_0=n\ZZ$.

\begin{theorem}
\label{thm:T_invertible}
$\Rep\left(T_q\right)$ is invertible in $\umor_2\left(\Mod_{\Rep\left(H\right)}\left(\Pre\right)\right)$.
\end{theorem}

\begin{proof}
The inclusion $\Vect_{\Lambda_0}\subseteq\Rep\left(T_q\right)$ is an $\EE_2$-map from
$\Rep\left(H\right)=\Vect_{\Lambda_0}$ with its symmetric monoidal structure, because
$Q|_{\Lambda_0}=1$, and it lands in transparent objects, because $\Lambda_0$ is the radical of $b$. So $\Rep\left(T_q\right)\simeq A^H$ for the fiber
$A\ceqq\Rep\left(T_q\right)\otimes_{\Rep\left(H\right)}\Vect$ (Proposition \ref{prop:equiv}(1)). The
fiber is computed on gradings, and is $\Vect_{\Lambda/\Lambda_0}$ with the induced braiding: a
pointed fusion category, which is non-degenerate because $\Lambda_0$ is the radical of $b$. As a
non-degenerate braided fusion category, it is invertible \cite[Theorem 1.1]{BJSS} (Remark \ref{rmk:abelian}). By Proposition \ref{prop:equiv}(2), $A^H$ is
invertible over $\Rep\left(H\right)$.
\end{proof}

\subsection{Outlook: the aggregated setting}
\label{sec:aggregated}

The algebraic setting cannot see analytic phenomena such as the Stone--von Neumann theorem, so it
is natural to ask which algebraic Witt classes become trivial once analytic objects are allowed. We
briefly describe one setting in which to ask this, and what can be said formally. In this section
$k=\CC$.

Aggregated vector spaces form a Grothendieck abelian category containing all complete bornological
spaces, e.g.\ Fr\'echet spaces, with a symmetric monoidal structure extending the completed
projective tensor product \cite{Steff:six, Schneiders}. Let $\Pre_{\Agg}$ be the $\infty$-category
of presentable $\infty$-categories enriched in its derived $\infty$-category, base changed to
$\CC$. A version $\QC^{\Agg}$ of $\QC$ on finitely generated affine derived $C^\infty$-schemes,
which include all compact manifolds and $\RR^n$, with a six-functor formalism satisfying K\"unneth
and base change, is constructed in \cite{Steff:six}; for compact $X$ it is aggregated modules over
$C^\infty\left(X\right)$. Base change $\Pre\to\Pre_{\Agg}$ is symmetric monoidal and preserves
colimits, so it induces a homomorphism of Witt groups (Definition \ref{defn:Witt})
\begin{equation}
\label{eqn:Witt_Agg}
\varphi\colon\mathrm{W}\left(\Pre\right)\longrightarrow\mathrm{W}\left(\Pre_{\Agg}\right) \ .
\end{equation}

For a compact torus $T$ (not the algebraic torus of Section \ref{sec:tori}), write
$\mathrm{Rep}\left(T\right)$ for the comodules over
$C^\infty\left(T\right)$ in aggregated modules, using
$C^\infty\left(T\right)\hat\otimes C^\infty\left(T\right)\simeq C^\infty\left(T\times T\right)$.

\begin{lem}
\label{lem:averaging}
$C^\infty\left(T\right)$ is descendable in $\mathrm{Rep}\left(T\right)$.
\end{lem}

\begin{proof}
Integrating against the normalized Haar measure defines a bounded map $C^\infty\left(T\right)\to\CC$,
and by translation invariance it is a map of comodules retracting the unit. So the unit is a retract
of $C^\infty\left(T\right)$, and in particular lies in the thick tensor ideal it generates
\cite{Mathew}.
\end{proof}

For a finite abelian group $\Gamma$ with a non-degenerate quadratic form $q$, we regard
$\Vect_\Gamma$ with the braiding determined by $q$ \cite{JS:braided}; it is invertible \cite{BJSS}
(Remark \ref{rmk:abelian}). Its central charge $c\in\ZZ/8$ is defined by
$\left|\Gamma\right|^{-1/2}\sum_{x\in\Gamma}q\left(x\right)=\exp\left(2\pi ic/8\right)$
\cite[Appendix 4]{MH}.

\begin{prop}
\label{prop:Witt_reduction}
Suppose given classes $a_V\in\mathrm{W}\left(\Pre_{\Agg}\right)$, one for each non-degenerate real
quadratic space $\left(V,b\right)$, such that
\begin{enumerate}
\item $a_V$ only depends on $\left(V,b\right)$ up to isometry;
\item $a_{V\dsum V'}=a_V+a_{V'}$;
\item for every even lattice $\left(N,b\right)$ with discriminant form $\left(\Gamma,q\right)$,
$\varphi\left[\Vect_\Gamma\right]+a_{N\otimes\RR}=0$.
\end{enumerate}
Then there is a class $g\in\mathrm{W}\left(\Pre_{\Agg}\right)$ with $8g=0$ such that
$\varphi\left[\Vect_\Gamma\right]=c\,g$ for every finite non-degenerate quadratic form
$\left(\Gamma,q\right)$ of central charge $c$. In particular, every such $\Vect_\Gamma$ with $c=0$
lies in the kernel of $\varphi$.
\end{prop}

\begin{proof}
Every $\left(\Gamma,q\right)$ is the discriminant form of an even lattice $\left(N,b\right)$
\cite[Theorem 1.10.1]{Nik}, so $\varphi\left[\Vect_\Gamma\right]=-a_{N\otimes\RR}$ by (3). Let
$h_\pm\ceqq a_V$ for $V=\RR$ with the form $\pm xy$. By Sylvester's law of inertia and (1), (2),
$a_V=r\,h_++s\,h_-$ for $V$ of signature $\left(r,s\right)$. The hyperbolic lattice is even with
trivial discriminant form, so (3) gives $h_++h_-=0$, and hence $a_V=\sigma\left(V\right)h_+$. The
$E_8$ lattice is even with trivial discriminant form and signature $8$, so $8h_+=0$. Finally
$\sigma\left(N\right)\equiv c$ modulo $8$ by Milgram's formula \cite[Appendix 4]{MH}, so
$g\ceqq-h_+$ works.
\end{proof}

So the relation $8g=0$ is forced by lattices, not by analysis. We expect such classes to arise as
$a_V=\left[\cA_V\right]$ for a braided category $\cA_V$ of sheaves on $V$ with braiding
$\exp\left(-2\pi i\,b\left(x,y\right)\right)$. Namely, for an even lattice $\left(N,b\right)$ and
$T=\left(N\otimes\RR\right)/N$, the form $b$ twists convolution on $\QC^{\Agg}\left(T\right)$ via a
biextension of $T\times T$, which has no algebraic counterpart since line bundles on algebraic tori
are trivial. Using Lemma \ref{lem:averaging}, the resulting $\EE_1$-algebra should Morita trivialize
its Drinfeld center $\cZ_N$, and a Peter--Weyl theorem for $T$, via Fourier transform, should
identify $\cZ_N$ with $\Vect_\Gamma\otimes\cA_{N\otimes\RR}$. This would give (3). Making this
precise requires defining $\cA_V$, including its decay conditions, and descent for module
categories along the descendable algebra of Lemma \ref{lem:averaging}.

\begin{question}
\label{qn:order_g}
Do such classes exist, and if so, what is the order of $g$? If $g=0$, which would be a $2$-shifted
analogue of the Stone--von Neumann theorem, every non-degenerate pointed braided fusion category lies in the kernel
of $\varphi$. Otherwise $\varphi$ detects part of the central charge modulo $8$.
\end{question}

\bibliographystyle{amsalpha}
\bibliography{references,missing_refs}
\end{document}